\documentclass[dvipsnames,letterpaper, 10 pt, journal]{ieeetran}

\IEEEoverridecommandlockouts	

\usepackage{algorithm}
\usepackage{xcolor}         
\usepackage{soul}
\usepackage{cite}
\usepackage{amsmath,amssymb,amsfonts,amsthm}
\usepackage{algpseudocode}
\usepackage{graphicx,color}
\usepackage{mathrsfs}
\usepackage{color}

  \newenvironment{smallarray}[1]
 {\null\,\vcenter\bgroup\small
  \arraycolsep=.13885em
  \hbox\bgroup$\array{@{}#1@{}}}
 {\endarray$\egroup\egroup\,\null}

\definecolor{myBlue}{RGB}{49,130,189}

\newtheorem{theorem}{Theorem}[section]
\newtheorem{lemma}[theorem]{Lemma}

\newtheorem{corollary}[theorem]{Corollary}

\newtheorem{remark}{Remark}
\newtheorem{assumption}[theorem]{Assumption}

\newcommand{\map}[3]{#1: #2 \rightarrow #3}

\usepackage{tabularx}
\usepackage{multirow}
\usepackage{makecell}

\newcommand\aamsout{\bgroup\markoverwith{\textcolor{violet}{\rule[0.5ex]{2pt}{1pt}}}\ULon}

\newcommand{\real}{\mathbb{R}}

\newcommand{\T}{\mathsf{T}} 

\newcommand{\mc}{\mathcal}

\newcommand{\trace}{\mathrm{tr}}

\newcommand{\expect}[1]{\mathbb{E}\left[#1\right]}
\newcommand{\Expect}[3]{\mathbb{E}_{#1}\left[#2 | #3 \right]}
\newcommand{\Expectover}[2]{\mathbb{E}_{#1}\left[#2\right]}
\newcommand{\Tr}[1]{\mathrm{tr}\left[#1\right]}

\newcommand{\overbar}[1]{\mkern 1.5mu\overline{\mkern-1.5mu#1\mkern-1.5mu}\mkern 1.5mu}

\newcommand{\1}{\mathds{1} }

\DeclareSymbolFont{bbold}{U}{bbold}{m}{n}
\DeclareSymbolFontAlphabet{\mathbbold}{bbold}

\newcommand\oprocendsymbol{\hbox{$\square$}}
\newcommand\oprocend{\relax\ifmmode\else\unskip\hfill\fi\oprocendsymbol}

\renewcommand{\theenumi}{(\roman{enumi})}

\newcommand*{\QEDA}{\hfill\ensuremath{\square}}%
\newenvironment{pfof}[1]{\vspace{1ex}\noindent{\itshape Proof of
    #1:}\hspace{0.5em}} {\hfill\QEDA\vspace{1ex}}

\graphicspath{{figs/}}
\makeatletter
\let\NAT@parse\undefined
\makeatother
\usepackage[colorlinks,urlcolor=blue, linkcolor=blue, citecolor=blue]{hyperref}

\title{Linear Dynamics meets Linear MDPs: Closed-Form Optimal Policies via Reinforcement Learning}
\author{Abed~AlRahman~Al~Makdah,~Oliver~Kosut,~Lalitha Sankar~and~Shaofeng Zou
  \thanks{A. A. Al Makdah, O. Kosut, L. Sankar, and S. Zou are with the School of Electrical, Computer and Energy Engineering at the Arizona State University,
    \href{mailto:aalmakda@asu.edu}{\{\texttt{aalmakda}},\href{mailto:okosut@asu.edu}{\texttt{okosut}},\href{mailto:lsankar@asu.edu}{\texttt{lsankar}},\href{mailto:zou@asu.edu}{\texttt{zou\}@asu.edu}}.}}

\usepackage{bibentry}

\begin{document}

\maketitle

\begin{abstract}
Many applications---including power systems, robotics, and economics---involve a dynamical system interacting with a stochastic and hard-to-model environment.  We adopt a reinforcement learning approach to control such systems. Specifically, we consider a deterministic, discrete-time, linear, time-invariant dynamical system coupled with a feature-based linear Markov process with an unknown transition kernel. The objective is to learn a control policy that optimizes a quadratic cost over the system state, the Markov process, and the control input.  Leveraging both components of the system, we derive an explicit parametric form for the optimal state-action value function and the corresponding optimal policy. Our model is distinct in combining aspects of both classical Linear Quadratic Regulator (LQR) and linear Markov decision process (MDP) frameworks. This combination retains the implementation simplicity of LQR, while allowing for sophisticated stochastic modeling afforded by linear MDPs, without estimating the transition probabilities, thereby enabling direct policy improvement.~We use tools from control theory to provide theoretical guarantees on the stability of the system under the learned policy and provide a sample complexity analysis for its convergence to the optimal policy. We illustrate our results via a numerical example that demonstrates the effectiveness of our approach in learning the optimal control policy under partially known stochastic dynamics. 
\end{abstract}

\section{Introduction}
In many applications, a well-modeled agent must interact with and make decisions in stochastic and hard-to-model environments, with the aim to optimize a certain cost that is affected by the agent's objective and the environment. A prominent example arises in power systems, where a controllable energy storage device evolves under known physical dynamics, yet must respond to uncertain net load demand driven by  exogenous factors such as variability in generation, consumer behavior, and weather conditions, all of which are unaffected by the device's control actions. Similar challenges appear for autonomous systems operating in unknown stochastic environments, or economic systems influenced by latent market factors. In such settings, designing optimal control strategies requires accounting for both the predictable evolution of the system and the stochastic nature of the surrounding environment. Effectively controlling such systems requires models that capture both the deterministic evolution of the agent’s state and the stochastic evolution of the environment. To address this challenge, in this work, we model the agent (e.g., battery energy storage system, self-driving car) with deterministic linear dynamics derived from first principles, while we model the environment (e.g., net load demand, traffic) as a linear Markov process.

Classical control theory offers elegant solutions for systems with entirely known dynamics, such as the Linear Quadratic Regulator (LQR) optimal control problem, which yields a closed-form optimal policy via Riccati equations \cite{BDA-JBM:07}. On the other hand, reinforcement learning (RL) approaches have developed data-driven techniques for decision-making in unknown environments, including model-based approaches. One compelling model in this setting is that of linear Markov decision processes (linear MDPs) that leverage feature-based representations to approximate the transition kernel \cite{SJB-AGB:96,FSF-MIR:07,RSS-AGB:18}. The linearity of the Markov kernel coupled with the non-linearity of feature functions results in a rich but tractable~model. 

Yet, even such a tractable framework  does not distinguish between the system dynamics and environment, viewing them as a single entity driven by the same dynamics. To this end, we propose an RL framework that combines both the LQR and linear MDP paradigms: a deterministic, discrete-time, linear time-invariant system coupled with a stochastic environment modeled as a feature-based linear Markov process that is not affected by the actions with unknown transition kernel. Our objective is to design a controller that optimizes a quadratic cost over the joint system and environment states and the control actions. In our setting, the deterministic part of the system and the quadratic cost weights are known a priori, only the environment’s transition kernel is unknown. By combining the structure of LQR and linear MDPs, we derive parametrized closed-form expressions for the optimal state-action value function and the corresponding optimal policy, capturing both deterministic dynamics and latent stochastic effects in a unified model. This hybrid model preserves the simplicity of LQR policies while incorporating the expressive stochastic modeling of linear MDPs. We use the least-squares value iteration (LSVI) algorithm to learn the parameters from online data in episodic fashion. The closed-form expression of the policy that optimizes the state-action value function makes it amenable to efficiently perform the policy update directly using the updated parameters at the end of each episode. Furthermore, because the unknown transition kernel is not affected by the control actions, our LSVI algorithm does not require exploration. We show that our LSVI achieves $\widetilde{\mathcal{O}}\left(T \sqrt{dL} \right)$ regret with high probability, where $d$, $T$, and $L$ denote the dimension of the feature-space of the linear Markov model, the time horizon of each episode, and the number of episodes, respectively. 
%
%
%
\subsection{Related work}
Our work lies at the intersection of optimal control and reinforcement learning, where we bridge ideas from the LQR optimal control problem and linear MDPs. Below, we review prior work that has been done in each area and highlight how our approach uniquely integrates them.

\paragraph{Linear Quadratic Regulator (LQR):} 
The classical LQR problem admits closed-form optimal control policies for linear systems. Traditional methods assume full knowledge of the system dynamics and cost, enabling the computation of optimal policies via Riccati equations \cite{BDA-JBM:07}. Recent work has studied the LQR problem in data-driven settings. Direct data driven approaches have been studied in \cite{CDP-PT:19,FD-PT-CDP:22,FC-GB-FP:23a}, where the optimal policy is learned directly from offline data generated by the open-loop system. Indirect data-driven approaches, explored in \cite{WA-DK-BJ-RM-MS:05,GRGDS-ASB-CL-LC:18,SD-HM-NM-BR-ST:20}, first identify a model of the system dynamics from data then solve the LQR problem using the identified model. Other works have studied the LQR problem in online learning setting \cite{MF-RG-SK-MM:18, HM-AZ-MS-MRJ:19,JB-AM-MF-MM:19,IF-BP:21,SJB-BEY-AGB:94}, where the optimal policy is learned online using policy gradient methods.


\paragraph{Reinforcement learning with function approximation:} In many reinforcement learning (RL) problems, the state or action spaces are too large (or continuous) to allow for tabular representations of value functions or policies \cite{JK-JAB-JP:13,VM-KK-DS-AG:13,SD-HA-CJ:2016}. To address this, function approximation techniques are employed to generalize from observed states and actions to unseen ones, enabling scalability and improved sample efficiency. Among the function approximation models, linear function approximation is particularly appealing due to its computational simplicity, theoretical tractability, and its ability to support efficient learning algorithms. Early approaches such as temporal difference learning, Q-learning, and least-squares temporal difference (LSTD) algorithms with linear value function approximation were explored in works like \cite{SJB-AGB:96,FSF-MIR:07,RSS-AGB:18}. While these methods laid important foundations, they often lacked sample efficiency guarantees and relied on heuristic exploration. Recent studies have introduced sample-efficient algorithms for linear MDPs, where the transition kernel is assumed to be a linear function of known features and unknown parameters \cite{LY-MW:19,LY-MW:20,CJ-ZY-ZW-MIJ:20}. In \cite{CJ-ZY-ZW-MIJ:20}, the authors developed a sample-efficient reinforcement learning algorithm for linear MDPs with a finite action space and a potentially infinite state space. Their model represents the transition kernel as a linear combination of known features with unknown probability measures, and assumes the reward function is linear in the same features with unknown parameters. In \cite{LY-MW:20}, the authors proposed a sample-efficient reinforcement learning algorithm under a linear MDP setting with possible infinite state and action space. In their framework, they assume the reward is known; further, their model introduces an additional structural assumption compared to \cite{CJ-ZY-ZW-MIJ:20}, by parameterizing the transition kernel with a low-dimensional unknown matrix. This assumption reduces the learning problem to estimating this matrix, thereby significantly lowering the overall learning complexity. In our framework, we model the stochastic environment as a feature-based linear Markov Process. Similar to \cite{CJ-ZY-ZW-MIJ:20}, we represent the transition kernel as a linear combination of known features with unknown probability measures. However, we assume a known quadratic cost that is independent of the features and consistent with the LQR framework, enabling efficient policy computation. This choice of the cost is more realistic and aligns with common formulations in engineering applications. Furthermore, our model avoids explicit parametric assumptions on the transition kernel made in \cite{LY-MW:20}, while allowing infinite state and action spaces. Additionally, our approach bypasses full model estimation by learning the value function directly through least-squares, benefiting from control-theoretic structure to ensure stability as well as computational and sample efficiency. Finally, our framework does not require exploration, since the environment is exogenous and is unaffected by the control inputs. 
Beyond linear MDPs, several works propose generalized model classes for sample-efficient RL including Bellman rank class, \cite{NJ-AK-AA-JL-RES:17}, linear Bellman-complete classes \cite{RM:05,AZ-AL-MK-EB:20}, witness rank  \cite{WS-NJ-AK-AA-JL:19}, and bilinear class \cite{SD-SK-JL-Sl-GM-WS-RW:21}.


\subsection{Contributions}

We list our contributions below.
\begin{itemize}
\item We propose an RL framework that unifies the classical LQR optimal control problem with linear MDPs. This hybrid model captures both deterministic dynamics of physical systems and stochastic evolution~of exogenous environments. To the best of our knowledge, this integration has not been addressed in existing literature. 

\item We derive a parametric form for the optimal state-action value function that decouples the agent’s dynamics from the environment’s stochasticity. This yields a closed-form policy that exhibits the LQR simplicity while inheriting the linear~MDPs' rich modeling capabilities. 

\item We propose a least-squares value iteration (LSVI) algorithm that learns the optimal policy by directly estimating the value function parameters. The LQR structure of our problem allows expressing the control policy in terms of the learned weights without having to optimize the value function at each time step as in \cite{CJ-ZY-ZW-MIJ:20}, thus simplifying the computational complexity of the algorithm.

\item We provide stability guarantees on the closed-loop system under the learned policy. These guarantees are given in terms of input-to-state stability, extending beyond standard sample-efficiency results found in RL literature, where it is typically assumed that the reward is bounded.

\item We derive a regret bound for our LSVI algorithm, which achieves a rate $\widetilde{\mathcal{O}}\left(T \sqrt{dL} \right)$ with high probability, where $d$, $T$, and $L$ denote the dimension of the feature-space of the linear Markov model, the time horizon of each episodes, and the number of episodes, respectively. 

\item We provide a numerical example to demonstrate the effectiveness of our framework, highlight its convergence, and verify the closed-loop stability of the learned~policy.

\end{itemize}

\section{Problem formulation}\label{sec: problem formulation}
Consider an agent obeying the discrete-time, linear, time-invariant dynamics over a finite time horizon
\begin{align}\label{eq: system}
 x_{t+1}=Ax_t+Bu_t , \quad t\in \{0,1,\ldots, T-1\},
\end{align}
where $x_t\!\in\!\mathcal{X}=\real^{n}$ denotes the state and $u_t \in \mathcal{U}= \real^{m}$ the input with $x_0 \sim \mathcal{N}(0,\Sigma_x)$ with $\Sigma_x \succ 0$. We assume the linear dynamical system, defined via the matrix pair $(A,B)$, is controllable.\footnote{When the system is controllable, it implies that there exist an input sequence, $u$, that can drive the system from its initial state, $x_0$ to any final state, $x_t$, within finite time horizon (see \cite[Section 9.8]{KO:10}).} We consider an environment evolving according to the discrete-time Markov~Process
%
\begin{align}\label{eq: environment}
 s_{t+1}| s_t\! \sim\! \mathbb{P}_t\left(s_{t+1}| s_t\right), \quad \!\! t\!\in\!\{0,\ldots,T\!-\!1\},
\end{align}
where $s\in \mathcal{S}\subset \real^p$ denotes the state of the Markov Process and $\mathbb{P}_t\left(s'|s\right)$ denotes the transition probability from state $s$ to $s'$, with $s_0 \sim \mu_0$ for some distribution $\mu_0\in \Delta(\mathcal{S})$, where $\Delta(\mathcal{S})$ denotes the set of distributions over $\mathcal{S}$. We assume that the matrices $A$ and $B$ in \eqref{eq: system} are known, while the transition probability, $\mathbb{P}_t$, in \eqref{eq: environment} is unknown. The agent follows a control policy $\map{\pi_t}{\mathcal{X}\times \mathcal{S}}{\mathcal{U}}$, where $u_t=\pi_t\left(x_t,s_t\right)$ is the action that the agent takes at state $x_t$ and $s_t$ at time $t$, for $t\geq0$. The objective is to find an optimal control policy, $\boldsymbol{\pi}=(\pi_0,\ldots,\pi_{T})$, that optimizes the following control task 
%
%
\begin{align}\label{eq: control task}
\begin{array}{ll}
   	\underset{\boldsymbol{\pi}}{\text{minimize}} & \displaystyle \expect{\sum_{t=0}^{T} c\left(x_t,s_t,u_t\right)}, \\
   	\text{subject to} & 
	x_{t+1} = Ax_t+Bu_t,\\
	& s_{t+1} \sim \mathbb{P}_t\left(s_{t+1}|s_t\right), \\
	 &u_t = \pi_t\left(x_t, s_t\right),
	  \end{array}
\end{align}
%
where $\map{c\!}{\mathcal{X}\!\times\!\mathcal{S}\! \times \mathcal{U}}{\real_{\geq0}}$ is the cost evaluated at $x_t$, $s_t$, and $u_t$ for $t\!\geq\! 0$ with $u_T=\pi_T(x_T,s_T)=0$. We restrict our search in \eqref{eq: control task} to the class of deterministic policies. We show later in Section \ref{subsec: value function} that the optimizer of \eqref{eq: control task} is indeed deterministic. We introduce the following assumptions on the transition probability in \eqref{eq: environment} and the cost~in~\eqref{eq: control task}.
%

%
\begin{assumption}{\bf \emph{(Linear Markov Process)}}\label{assump: linMDP}
Let $\map{\phi}{\!\mathcal{S}}{\!\real^d}$ be a known feature vector and $\mu_t\!\in\! \mathbb{R}^d$ a vector of $d$ unknown signed measures over $\mathcal{S}$. For $s'\!,s\!\in\! \mathcal{S}$, we have
\begin{align}\label{eq: linApprox}
 \mathbb{P}_t\left(s'|s\right)=\phi(s)^{\T} \mu_t(s'),
\end{align}
We assume $\|\phi(s)\|\leq 1/\sqrt{d}$ and $\|s\|\leq \delta_s$ for all $s \in \mc{S}$, and $\expect{\phi(s_t)\phi(s_t)^{\T}}\succ 0$ for all $t$.
\end{assumption}

\begin{assumption}{\bf \emph{(Quadratic cost)}}\label{assump: QuadCost}
For $x\in\mathcal{X}$, $s\in \mathcal{S}$, and $u \in \mathcal{U}$, we have
\begin{align}\label{eq: QuadCost}
c(x,s,u)=
 \begin{bmatrix}
 x \\  s \\ u
\end{bmatrix}^{\T}
\underbrace{\begin{bmatrix}
W & F & D \\
 F^{\T} & M & H \\
 D^{\T}    & H^{\T}    & R \\
\end{bmatrix}}_{P}
 \begin{bmatrix}
 x \\  s \\ u
\end{bmatrix}, 
\end{align}
where $P\succeq 0$ is known and $R\succ 0$. Further, we assume the pair $(A,W^{1/2})$ is observable.
\end{assumption}

Assumption \ref{assump: linMDP} is inspired by the linear MDP framework introduced in \cite{SJB-AGB:96,FSF-MIR:07,CJ-ZY-ZW-MIJ:20}. However, unlike the original definition, our model assumes that the stochastic process governs only the exogenous state and is unaffected by control input. This assumption is motivated by the fact that, in our target applications, the environment is not influenced by control actions. Moreover, it simplifies the expression of the optimal policy, as the optimal policy requires minimizing a quadratic function in the input $u$ (from Assumption \ref{assump: QuadCost}), rather than the nonlinear (possibly non-convex) function $\phi$.
%

We aim to learn the optimal policy for the system in \eqref{eq: system} that optimizes the control objective defined in Problem \eqref{eq: control task}. This objective is also affected by a stochastic process, $s_t$, which evolves according to the linear model defined in \eqref{eq: linApprox}. The system dynamics matrices, $A$ and $B$, in \eqref{eq: system}  as well as the weight matrices and the feature map in \eqref{eq: linApprox} are assumed to be known, while the probability measures, $\mu$ in \eqref{eq: linApprox} are assumed to be unknown. 
%
%
%
%
We define the value function $\map{V_{t}^{\pi}}{\mathcal{X}\! \times\! \mathcal{S}}{\!\real}$ as the expected cumulative cost incurred~under policy $\pi$ starting from state $x_t$ and $s_t$ at time $t\!\geq\! 0$, given by
%
\begin{align*}
 V_{t}^{\pi} \left(x,s\right) \triangleq  \expect{\sum_{i=t}^{T} c \left(x_{i},s_{i}, \pi_i\left(x_{i},s_{i}\right)\right) \Big| x_t = x, s_t = s}.
\end{align*}
%
Further, we define the state-action value function $\map{Q_{t}^{\pi}}{\mathcal{X}\!\times\! \mathcal{S}\!\times\! \mathcal{U}}{\real}$~as the expected cumulative cost under policy $\pi$ starting from state $x_t$, $s_t$, and action $u_t$ at time $t\!\geq\! 0$, given~by
%
\begin{align*}
Q_{t}^{\pi}&\!\left(x,s,u\right)\triangleq c\! \left(x,s, u \right)\\
&+\expect{\sum_{i=t+1}^{T} \!c\! \left(x_{i},s_{i}, \pi_t\left(x_{i},s_{i}\right)\right) \Big| x_t\!=\!x, s_t\!=\!s , u_t\!=\!u}.
\end{align*}
To learn the optimal policy, we focus on estimating the state-action value function $Q_{t}^{\pi}$, since it directly guides policy improvement through greedy action selection. In particular, by learning~an appropriate parametric approximation of the state-action value function, $Q_{t}^{\pi}$, we can infer an optimal policy without explicitly learning the transition probability measures, $\mu$, in \eqref{eq: linApprox}. This approach leverages the structure of the system and cost, allowing us to bypass the need for full system identification and instead focus on value function approximation within the RL framework.
\begin{remark}{\bf \emph{(On the knowledge of $A$ and $B$ in \eqref{eq: system})}}\label{rmrk: A and B}
In many control applications—such as robotics, aerospace, and power systems—the physical plant dynamics (i.e., $A$ and $B$) can be easily derived from first principles or can be accurately identified through standard system identification techniques prior to deployment. Our framework leverages this knowledge to focus on learning the stochastic environment component, which simplifies the computational complexity of the policy update step, and enables stability-aware control without requiring aggressive exploration. Nonetheless, our approach can be extended to settings where $A$ and $B$ are unknown, which we leave for future~work. In addition, incorporating additive process noise into the system dynamics is another natural extension that we plan to explore.
\end{remark}

\section{Main results}
We now present our main results. First, we leverage the linear structures of the system~in~\eqref{eq: system}, the transition model in \eqref{eq: linApprox}, and the quadratic structure of the cost in \eqref{eq: QuadCost} to derive a parametric expression for the state-action value function that is linear in the feature map, $\phi$, along with a parametric expression for the corresponding optimal greedy policy. Second, we introduce a least-squares value iteration algorithm to learn the parameters of the state-action value function, and therefore learn the optimal policy. Third, we provide stability guarantees on the closed-loop system under the learned policy. Finally, we provide convergence analysis for our algorithm and derive an upper bound on the regret.

\subsection{State-action value function approximation}\label{subsec: value function}
Let the optimal value function at time $t$ and evaluated at $x\in \mc{X}$ and $s\in \mc{S}$ under the optimal policy, $\pi_t^*$, be denoted by $V_t^*(x,s)$. Following the Bellman optimality equation, we can write the optimal state-action value function at time $t$ and evaluated at $x\in \mc{X}$, $s\in \mc{S}$, and $u\in \mc{U}$ under $\pi_t^*$ as 
\begin{align*}
 Q_t^*(x,s,u)\!=\!c(x,s,u) \!+\! \mathbb{E}_{s'\sim \mathbb{P}_t(s'|s)\!\!}\left\{ V_{t+1}^*(Ax\!+\!Bu,s')| s\right\}\!.
\end{align*}
%

The next result provides an explicit parametric form for the state-action value function $Q_t$.
\begin{theorem}{\bf \emph{($Q$-function representation)}}\label{thrm: Qfunction}
Consider the dynamics in \eqref{eq: system} and the Markov Process in \eqref{eq: environment}. Let Assumption \ref{assump: linMDP} and Assumption \ref{assump: QuadCost} be satisfied. Then, for any $x\in \mc{X}$, $s\in \mc{S}$, and $u\in \mc{U}$, and under $\pi_t^*$ for $t\geq0$, there exists $\overline{h}_{i,t+1}\in \mathbb{R}^{n}$ and $\overline{q}_{i,t+1}\in\mathbb{R}$ such that
\begin{align} \label{eq: Qfunc} 
\begin{split}
 Q_t^*(x,s,u)=& c(x,s,u)\!+\! \left(Ax+Bu\right)^{\T} G_{t+1}\left(Ax+Bu\right)\\
 &+\sum_{i=1}^{d}\phi_i(s)\Big(2\left(Ax\!+\!Bu\right)^{\T} \overline{h}_{i,t+1}+\overline{q}_{i,t+1}\Big),
 \end{split}
\end{align}
where $G$ solves the discrete-time algebraic Riccati equation
\begin{align}
\begin{split}
 G_t=&A^{\T}G_{t+1}A + W \\
 - (A^{\T}&G_{t+1}B \!+\! D) (R\!+\!B^{\T}G_{t+1}B )^{-1} (B^{\T}G_{t+1}A \!+\! D^{\T}), \label{eq: riccati_eq}
 \end{split}
 \end{align}
 with $G_T=W$.
\end{theorem}

A proof of Theorem \ref{thrm: Qfunction} is in Appendix \ref{app: A}. Several comments are in order. First, by leveraging the linearity of the system in \eqref{eq: system} and the Markov process in \eqref{eq: linApprox}, along with the quadratic structure of the cost in \eqref{eq: QuadCost}, the state-action value function in \eqref{eq: Qfunc} exhibits a structure that decouples the linear system state $x$ and action $u$ from the exogenous state $s$. Second, the derived expression of the state-action value function in \eqref{eq: Qfunc} is linear in the feature map $\phi$ and the weight parameters $\overline{h}_{i,t}$ and $\overline{q}_{i,t}$. Third, the weight parameters $\overline{h}_{i,t}$ and $\overline{q}_{i,t}$ depend on the unknown transition probability $\mathbb{P}(\cdot|s)$ in \eqref{eq: linApprox}, and therefore, learning the state-action value function boils down to learning these weights, thereby bypassing the need to explicitly learn the probability measures, $\mu$, in \eqref{eq: linApprox}.

The optimal policy is found by minimizing the $Q$ function over the input $u$. Since by Theorem~\ref{thrm: Qfunction}, $Q$ is quadratic in $u$, this optimal policy can be found in closed form, as shown in the following corollary, which expresses the optimal policy in terms of feedback gains and weight parameters.

\begin{corollary}{\bf \emph{(Optimal policy representation)}}\label{cor: opt_policy}
For any $x\in \mc{X}$, $s\in \mc{S}$, and~$t \in \{0,1,\cdots, T-1\}$
\begin{align*}
\begin{split}
u_t^*(x,s)=&\pi_t^*(x,s)\\
=&K_{t,x}x + K_{t,s}s + K_{t,h}\sum_{i=1}^{d} \phi_i(s) \overline{h}_{i,t+1},
\end{split}
\end{align*}
where $\overline{h}_{i,t+1}$ is as in Theorem \ref{thrm: Qfunction} and
\begin{align}\label{eq: feedback gains}
\begin{split}
 K_{t,x}&=-\left(R +B^{\T} G_{t+1} B \right)^{-1} \left( B^{\T}G_{t+1}A+D^{\T}\right),\\
 K_{t,s}&=-\left(R +B^{\T} G_{t+1} B \right)^{-1}H^{\T},\\
 K_{t,h}&=-\left(R +B^{\T} G_{t+1} B \right)^{-1}B^{\T},
 \end{split}
\end{align}
and $G_{t+1}$ satisfies \eqref{eq: riccati_eq}.
\end{corollary}

\begin{algorithm} 
    \caption{Least-Squares Value Iteration}
    \label{alg: alg1}
    \begin{algorithmic}[1]
     \State Given: $L$, $R_{\theta}$, $\lambda$
            \For{episode $\ell=1,\cdots , L$\hspace{5pt}}\\
                \hspace{12pt} $x^{\ell}(0) \overset{\text{i.i.d}}{\sim} \mathcal{N}(0,\Sigma_x)$ with $\Sigma_x \succ 0$\\
                \hspace{12pt} $s^{\ell}(0) \overset{\text{i.i.d}}{\sim} \mu_0$ such that $\expect{\phi(s_0)\phi(s_0)^{\T}}\succ 0$ 
                \For{step $t=T-1,\cdots , 0$\hspace{5pt}}
                    \State {\small{$\Lambda_t^{\ell}\!\gets\! {\!\sum_{i\!=\!1}^{\ell\!-\!1}\!\!{Y\!(x^{i}_t,\!u^{i}_t)}^{\T}\! \phi(s^{i}_t) {\phi(s^{i}_t)}^{\T}\! Y\!(x^{i}_t,\!u^{i}_t) \! +\! \lambda I_{dn+d}}$}}
                    \State {\small{$\theta_{t+1}^{\ell}\! \gets\! (\Lambda_t^{\ell})^{-1} \!\!\sum_{i\!=\!1}^{\ell\!-\!1}\!{Y\!(x^{i}_t,\!u^{i}_t)}^{\T}\! \phi(s^{i}_t)\epsilon_{t+1}^{\ell}\!(x^{i}_{t+1},\!s^{i}_{t+1})$}}
                    \If{$\|\theta_{t+1}^{\ell}\|>R_{\theta}$}
			\State $\theta_{t+1}^{\ell} \gets \frac{R_{\theta}}{\|\theta_{t+1}^{\ell}\|}\theta_{t+1}^{\ell}$
		   \EndIf
                \EndFor
                
                \For{step $t=0,\cdots , T-1$\hspace{5pt}}
                \State $u^{\ell}_t \gets K_{t,x} x^{\ell}_t + K_{t,s} s^{\ell}_t + K_{t,h} ({\phi(s^{\ell}_t)}^{\T} \otimes Z )\theta_{t+1}^{\ell}$
                \State Take action $u^{\ell}_t$
                \State Observe $x^{\ell}_{t+1}$ and $s^{\ell}_{t+1}$
                \EndFor
            \EndFor
    \end{algorithmic}
\end{algorithm}

\subsection{Learning weight parameters of the value function via least-squares value iteration}
In this subsection, we learn the weight parameters, $\overline{h}$ and $\overline{q}$, that parameterize the state-action value function in Theorem \ref{thrm: Qfunction}. To this aim, we propose a least-squares value iteration algorithm (Algorithm \ref{alg: alg1}) that is inspired by \cite{CJ-ZY-ZW-MIJ:20}. Before we lay out the steps of our algorithm, we introduce the following notations. At each time step, $t$, we concatenate the parameters $\overline{h}_{i,t}$ and $\overline{q}_{i,t}$ for $i\in\{1,\cdots,d\}$ as 
\begin{align}\label{eq: theta}
\theta_t= \begin{bmatrix}
 \theta_{1,t}^{\T}&
 \cdots&
\theta_{d,t}^{\T}
\end{bmatrix}^{\T},
\quad
\text{where}
\quad
\theta_{i,t} =\begin{bmatrix}
 \overline{h}_{i,t}\\
 \overline{q}_{i,t}
\end{bmatrix}.
\end{align}
Using the notation in \eqref{eq: theta}, we re-write the value function in Theorem \ref{thrm: Qfunction} and the policy in Corollary \ref{cor: opt_policy}, respectively,~as
\begin{align}\label{eq: policy alg}
 \begin{split}
 Q_t(x,s,u)=& c(x,s,u)\!+\! \left(Ax+Bu\right)^{\T} G_{t+1}\left(Ax+Bu\right)\\
 &+ {\phi(s)}^{\T} Y(x,u) {\theta}_{t+1}, \\
 u_t(x,s)=&K_{t,x}x + K_{t,s}s + K_{t,h} \left({\phi(s)}^{\T} \otimes Z\right) \theta_{t+1},
 \end{split}
\end{align}
where $Y(x,u)\!=\!I_d\!\otimes\! [2 \left(Ax + Bu\right)^{\T} , 1]$ and $Z\!=\![I_n, 0_{n \times 1}]$. Now we lay out the steps of our least-squares value iteration algorithm (Alg. \ref{alg: alg1}). Our algorithm consists of an outer loop over $L$ episodes, where each episode consists of two loops: 1) backward-in-time weight update loop (lines 5-11) and 2) forward roll-out and data collection loop (lines 12-16). During the first pass of episode~$\ell$ (lines 5–11), we treat the data collected in the previous $\ell-1$ episodes as a fixed dataset
\begin{align}\label{eq: dataset}
 \mathcal D_{\ell-1}\!:=\!\bigl\{
        (x^{i}_t,s^{i}_t,u^{i}_t,x^{i}_{t\!+\!1},s^{i}_{t\!+\!1})
        :i<\ell,\;0\leq t<T
     \bigr\}.
\end{align}
At each time step $t$, $\theta$ minimizes a regularized least-squares loss---the squared error between the parametric state-action value function in \eqref{eq: policy alg} and the Bellman target (immediate cost plus the estimated value of the next state). Solving this problem on past trajectory data yields an accurate value-function approximation and enables closed-form greedy policy updates without estimating the transition probabilities. The regularized least-squares regression is stated as
 \begin{align*}
\begin{split}
  \theta^{\ell}_{t+1}\!\!=\!&
      \underset{\theta\in\mathbb R^{d(n+1)}}{\arg\min}
      \sum_{i=1}^{\ell-1}
            \!\Bigl(
                \phi(s_t^{i})^{\T} Y(x_t^{i},u_t^{i}) \theta -\!
               \epsilon_{t+1}^{\ell}(x^i_{t+1},s^i_{t+1})
            \Bigr)^{2}\\
      &\qquad \qquad\qquad+\lambda\|\theta\|^{2},
      \end{split}
\end{align*}
where
\begin{align*}
\epsilon_{t+1}^{\ell}\left(x,s\right) =& 2\left(x\right)^{\T}h^{\ell}_{t+1}\left(s\right) + q^{\ell}_{t+1}\left(s\right),\\
h^{\ell}_{t+1}\left(s\right)=&\left(A^{\T} + K_{t,x}^{\T}B^{\T}\right)(\phi(s)^{\T} \otimes Z) \theta^{\ell}_{t+2}\\
&+ \left(F +K_{t,x}^{\T}H^{\T}\right)s,\\
q^{\ell}_{t+1}(s)=&\left({\phi(s)}^{\T} \otimes \overline{Z}\right) \theta^{\ell}_{t+2} +s^{\T}(M + H K_{t,s}) s\\
 &+ {\theta^{\ell}_{t+2}}^{\T} \left({\phi(s)}\!\otimes\! Z^{\T}\right) B K_{t,h} ({\phi(s)}^{\T} \!\otimes\! Z) \theta^{\ell}_{t+2}\\
 &+2 s^{\T}H K_{t,h}\left({\phi(s)}^{\T} \otimes Z\right) \theta^{\ell}_{t+2},
\end{align*}
with $\overline{Z}=[0_{1\times n}, 1]$. Unlike prior work (e.g., \cite{CJ-ZY-ZW-MIJ:20}), we leverage the structure of our model to derive a closed-form expression for the Bellman target in terms of previously learned parameters, thereby avoiding an inner optimization over the action space at each time step (often required in discrete action space settings). In fact, $\epsilon_{t+1}^{\ell}\left(x,s\right)$ is obtained directly from this closed-form Bellman target (see Appendix \ref{app: LSVI}). The closed-form parameter update is given by
\begin{align}\label{eq: LSVI sol}
\begin{split}
      \Lambda^{\ell}_{t}=&
     \sum_{i=1}^{\ell-1}
        Y(x^{i}_t,u^{i}_t)^{\T}\phi(s^{i}_t)\phi(s^{i}_t)^{\T}Y(x^{i}_t,u^{i}_t) +\lambda I_{d(n+1)}
        ,\\
     \theta^{\ell}_{t+1}=&\bigl(\Lambda^{\ell}_{t}\bigr)^{-1}\sum_{i=1}^{\ell-1}{Y\left(x^{i}_t,u^{i}_t\right)}^{\T} \phi\left(s^{i}_t\right)\epsilon_{t+1}^{\ell}\left(x^{i}_{t+1},s^{i}_{t+1}\right),
     \end{split}
\end{align}
which recover lines 6 and 7 of Alg. \ref{alg: alg1}. For $\ell=1$, we set $\theta^{\ell}_{t+1}=0$ and $\Lambda_t^{\ell}=\lambda I_{d(n+1)}$ for $t\in \{0,\cdots, T-1\}$. The regularizer term $\lambda I_{d(n+1)}$ ensures numerical stability, the projection step in lines 8-10 makes sure that the norm of the learned parameters is uniformly bounded for $t\in \{0,\cdots,T-1\}$ and $\ell \in\{1,\cdots ,L\}$. In the second pass (lines 12–16) the newly computed parameters $\theta^{\ell}_{t+1}$ are plugged into the greedy closed-form policy \eqref{eq: policy alg},  
\begin{align*}
      u^{\ell}_t(x_t^{\ell},s_t^{\ell})\!=\!K_{t,x}x^{\ell}_t\!+\!K_{t,s}s^{\ell}_t
                 +K_{t,h}\bigl(\phi(s^{\ell}_t)^{\!\top}\!\otimes \!Z\bigr)\theta^{\ell}_{t+1},
\end{align*}
generating a new trajectory $(\{x^{\ell}_t,s^{\ell}_t,u^{\ell}_t\}_{t=0}^{T})$. These samples are appended to the collected data \eqref{eq: dataset}, and will be used in the next episode’s backward update. Notice that Alg. \ref{alg: alg1} does not require exploration as in \cite{CJ-ZY-ZW-MIJ:20}, which we discuss in the following remark.


%
\begin{remark}{\bf \emph{(Role of exploration)}}\label{rmrk: exploration}
In classical reinforcement learning setting, exploration (e.g., using $\varepsilon$-greedy or optimism-based methods) is necessary to sufficiently explore the environment and estimate unknown transition dynamics. However, our framework does not require exploration. This is because the stochastic component of the environment is modeled as an exogenous Markov process that evolves independently of the control inputs (see Assumption \ref{assump: linMDP}), and the system dynamics, $A$ and $B$, are known. Our algorithm estimates the value function parameters, $\theta$, via a least-squares procedure using observed trajectories without the need to infer the transition probabilities explicitly. As a result, the optimal policy can be computed in closed form by minimizing a known quadratic function of the control input.
\end{remark}
\begin{remark}{\bf \emph{(Choice of $R_{\theta}$)}}\label{rmrk: R_theta}
The projection radius $R_{\theta}$ in Alg. \ref{alg: alg1} ensures that the learned parameters at each episode and time step remain within a ball of radius $R_{\theta}$, which ensures numerical stability. Moreover, it plays a crucial role in the theoretical analysis (i.e., stability and regret bound). In practice, $R_{\theta}$ should be chosen large enough to contain the true parameters, $\theta^*$, but not excessively large in order to keep the constants in the stability and regret bounds moderate. In Appendix~\ref{app: true theta}, we derive an upper bound on $\|\theta^*\|$; if $R_{\theta}$ is chosen to be larger than this bound, then the ball of radius $R_{\theta}$ is guaranteed to contain $\theta^*$. In particular, we show that $\theta^*$ is contained in this ball if $R_{\theta}\! \geq\! c_{\theta}\sqrt{d}$, where $c_{\theta}\!>\!0$ depends on known problem parameters, e.g., the system matrices, cost weights, the feature map, and the bound~on~$s$.
\end{remark}
%
%
%
%
%
%
%
%
%
%
%

%
\subsection{Input-to-State Stability}
It is critical to ensure that the learned policy stabilizes the closed-loop system in each episode, particularly in settings where the environment evolves independently of the control actions and safety is a concern. To this end, we establish an input-to-state stability (ISS) bound for the system under the learned policy, which we present in the following result.
%
\begin{theorem}{\bf \emph{(Input-to-state stability)}}\label{thrm: ISS}
Consider system \eqref{eq: system}, let $u$ be the output of Algorithm \ref{alg: alg1} at episode $\ell$. Let $\|\theta_t^{\ell}\|\leq R_{\theta}$, $\|K_{t,s}\|\leq \overline{K}_s$, and $\|K_{t,h}\|\leq \overline{K}_h$ for $t\in\{0, \cdots , T - 1\}$. Let $x^{\ell}_0$ be the initial state in episode $\ell$. Then, under Assumptions \ref{assump: linMDP}~and~\ref{assump: QuadCost},
%
\begin{align}\label{eq: ISS}
  ||x^{\ell}_t||\leq  \alpha\rho^t \|x^{\ell}_0\| + \frac{\alpha\|B\|}{1-\rho}\left(\overline{K}_s\delta_s + \frac{\overline{K}_hR_{\theta}}{\sqrt{d}}\right),
\end{align}
for $t\!\in\!\!\{0,\!\cdots\!, T\!-\!1\}$, where $\alpha\!>\!0$ and $0\!<\!\rho\!<\!1$ are~constants. 
\end{theorem}
A proof of Theorem \ref{thrm: ISS} is deferred to Appendix \ref{app: ISS}. Several comments are in order. First, Theorem \ref{thrm: ISS} implies that the state trajectory at each episode $\ell \in \{1,\!\cdots\! , L\}$ remains bounded in terms of the initial condition, linear system dynamics, the control gains in Corollary \ref{cor: opt_policy}, and $R_{\theta}$. Second, the first term on the right-hand side of \eqref{eq: ISS}, which depends on the initial state, decays exponentially with time, while the second term is independent of time and the number of episodes. This latter term depends on the system matrices, feedback gains, the bound on $s$, $R_{\theta}$ from Algorithm \ref{alg: alg1}, and the dimension of the feature map, $d$. This result leverages the known system dynamics and the structure of the control policy, extending traditional stability notions in control to learning-based policies in partially known environments.

\subsection{Regret analysis}\label{sec: main results}
We analyze the convergence of Algorithm \ref{alg: alg1} by bounding the regret accumulated over multiple episodes. In particular, we define the regret $\mc{R}(L)$ as the difference between the total cost incurred by the learned policy and that of the optimal policy over $L$ episodes. Mathematically, for $L$ episodes, the total expected regret is defined as
\begin{align}\label{eq: regret}
\mc{R}(L)=\sum_{\ell=1}^{L}\left(V_0^{\ell}(x_0^{\ell},s_0^{\ell})-V_0^*(x_0^{\ell},s_0^{\ell})\right).
\end{align}
where $V_0^{\ell}(x_0^{\ell},s_0^{\ell})$ denotes the value evaluated at the initial states $x_0^{\ell}$ and $s_0^{\ell}$ under the policy learned at episode $\ell$, and $V_0^*(x_0^{\ell},s_0^{\ell})$ is the value of the optimal policy evaluated at the initial states $x_0^{\ell}$ and $s_0^{\ell}$. We derive a bound on the regret in the following result.
%
%
%
%
\begin{theorem}{\bf \emph{(Regret bound)}}\label{thrm: regret bound}
Let Assumptions \ref{assump: linMDP} and~\ref{assump: QuadCost} be satisfied. Let $\|Y(x_t^{\ell},u_t^{\ell})^{\T}\phi(s^{\ell}_t)\|\leq \delta_{\psi}$, for $t\in\{0,\cdots,T\}$ and $\ell \in \{1,\cdots,L\}$. Let $\beta=\log\left(1+\frac{L\delta_{\psi}^2}{\lambda}\right)$ with $\lambda>0$. Let $\delta \in [0,1/3]$. Then, with probability at least $1-3\delta$
\begin{align*}
&\mathcal{R}(L) \leq \sigma \sqrt{2TL\log(1/\delta)}\\
&\!+\!\delta_{\psi}T(\frac{1}{\sqrt{\lambda}}\!+\!\frac{4\sqrt{L}}{\sqrt{\gamma}}) \Big(\!\sigma\sqrt{2dn\beta \!+\! 2\log(\frac{1}{\delta})}\!+\!(R_{\theta}\!+\!2\delta_{v})\sqrt{\lambda}\Big)
\end{align*}
%
%
%
where $\sigma>0$, $\gamma>0$ and $\delta_{v}>0$ are constants that do not depend on $L$ and $T$, and do not scale with $d$. Further, $\delta_{\psi}$ scales with $\mathcal{O}(1/\sqrt{d})$ and $R_{\theta}$ scales with $\mathcal{O}(\sqrt{d})$.
\end{theorem}

A proof of Theorem \ref{thrm: regret bound} is presented in Appendix \ref{app: regret bound}. Theorem \ref{thrm: regret bound} provides a probabilistic upper bound on the cumulative regret. Several comments are in order. First, the leading term of the bound scales as $\mathcal{O}\left(T\sqrt{dL\log(L)}\right)$ or $\widetilde{\mathcal{O}}\left(T \sqrt{dL} \right)$, which matches, in terms of the number of episodes, the rate reported in \cite{CJ-ZY-ZW-MIJ:20}. Second, our bound grows linearly in $T$ and $\sqrt{d}$, in contrast to the $T^2$ and $d\sqrt{d}$ factors in \cite{CJ-ZY-ZW-MIJ:20}, respectively.
Third, the constants $\sigma$ and $\delta_v$ are independent of $L$ and $T$, and they do not scale with $d$ and they depend only on the system matrices, the cost in \eqref{eq: system}, the cost weight matrices in \eqref{eq: QuadCost}, the bound on the state, $x_t$, in \eqref{eq: ISS}, and the bound on the exogenous state $s_t$ in \eqref{eq: linApprox}. In fact, they arise from the uniform upper bound on the value function, $V_t$ (see Appendix~\ref{app: regret bound} for the details). Fourth, the constant $\gamma$ satisfies $\expect{Y(x_0,u_0)^{\T}\phi(s_0)\phi(s_0)^{\T}Y(x_0,u_0)}\!\succeq\! \gamma I_{d(n+1)}$,~which holds because the initial states, $x_0$ and $s_0$ are drawn independently in each episode. Finally, the bound in Theorem~\ref{thrm: regret bound} suggests, when the initial states, $x_0$ and $s_0$, are fixed~for~all episodes, Alg. \ref{alg: alg1} can learn an $\varepsilon-$optimal policy, $\pi$, that satisfies $V_0^{\pi}(x_0,\!s_0)\!-\!V_0^{*}(x_0,\!s_0)\!\leq\!\varepsilon$ after~$L\!=\!\widetilde{\mathcal{O}}\!\left(\frac{dT^2}{\varepsilon^2} \right)$~episodes.
\begin{figure*}[!t]
  \centering
  \includegraphics[width=1\textwidth,trim={0cm 0cm 0cm
    0cm},clip]{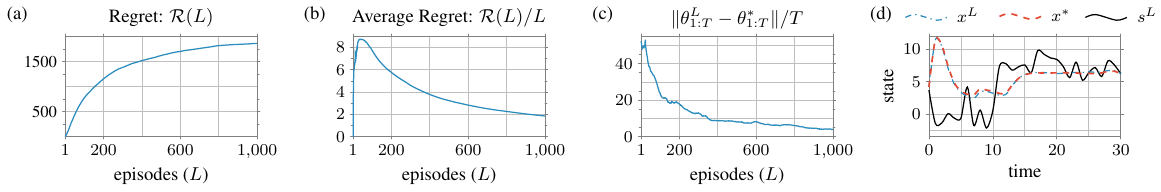}
  \caption{This figure shows the numerical results for the setting described in Section~\ref{sec: example}. Panel (a) shows the regret as a function of $L$. We observe that the regret scales as $\widetilde{\mathcal{O}}(\sqrt{L})$, aligning with Theorem~\ref{thrm: regret bound}. Panel (b) shows the average regret as a function of $L$, and we observe that it converges as  $L$ increases. Panel (c) shows the norm of the estimation error between the learned and the true parameters, averaged over the episode horizon $T$, as a function of $L$. We observe that the estimation error decreases with $L$, indicating that the learned policy converges to the optimal one. Panel (d) shows the state trajectory generated by the system in \eqref{eq: ex_dynamics} under the learned policy at episode $L=1000$ (dot-dashed blue line) and under the optimal policy (dashed red line). It also shows the exogenous state trajectory generated by the linear Markov process in \eqref{eq: ex_MDP} (solid black line). We observe that the trajectory under the learned policy closely matches that of the optimal policy, and both track the mean of the exogenous~state.
}
    \label{fig: example}
\end{figure*}
\section{Numerical Example} \label{sec: example}
We consider a discrete-time, linear, time-invariant system
\begin{align}\label{eq: ex_dynamics}
 x_{t+1}=
\underbrace{ \begin{bmatrix}
 1.8 & 1.2\\
 0 & 1.19
\end{bmatrix}}_{A}x_t
+
\underbrace{\begin{bmatrix}
 0\\1
\end{bmatrix}}_{B}u_t,
\end{align}
and a stochastic environment evolving according to a feature-based linear Markov process with
\begin{align}\label{eq: ex_MDP}
\begin{split}
s_{t+1}|s_t\! \sim\!\underbrace{ 
\begin{bmatrix}
\frac{f_1(s_t)}{f_1(s_t)+f_1(s_t)} & \frac{f_2(s_t)}{f_1(s_t)+f_1(s_t)}
\end{bmatrix}}_{{\phi(s_t)}^{T}}\!\!
\underbrace{\begin{bmatrix}
 \mc{N}(s_{t+1};7,1)\\
  \mc{N}(s_{t+1};-1,1.5)
\end{bmatrix}}_{\mu_{t+1}(s_{t+1})}\!,
\end{split}
\end{align}
where $f_1(s_t)\!=\!\exp\left(\frac{-\left(s_t-\nu_1\right)^2}{2\rho_1^2}\right)$ and $f_2(s_t)=\exp\left(\frac{-\left(s_t-\nu_2\right)^2}{2\rho_2^2}\right)$, where $\nu_1=7$, $\nu_2=-1$, $\rho_1=5$, and~$\rho_2=3$, with $\|s_t\|\leq \delta_s=15$ for all $t$. We define the cost function to capture the tracking error between the first state of \eqref{eq: ex_dynamics} and the exogenous state, $s$, and~is expressed~as
\begin{align*}
 c(x,s,u)=&\left(Cx-s\right)^{\T} M \left(Cx-s\right) +u^{\T}Ru\\
 =&x^{\T}\underbrace{C^{\T}MC}_{W}x + s^{\T} M s + u^{\T}R u -2s\underbrace{MC}_{F^{\T}}x,
\end{align*}
where $C=[1, 0]$, $M=1$, and $R=1$. First, we use the matrices $A$ and $B$, along with the cost weight matrices, to compute the feedback gains in Corollary~\ref{cor: opt_policy}. Then, we apply Algorithm~\ref{alg: alg1} to learn the parameters, $\theta$, using $L=1000$ episodes, each with horizon $T=30$. We set $\lambda=2$ and $R_{\theta}=500$ (see Remark \ref{rmrk: R_theta}). At each episode, we sample $x_0\overset{i.i.d.}{\sim} \mathcal{N}(3,1)$ and $s_0\overset{i.i.d.}{\sim} \mathcal{N}(3,1)$, which are independent of each other. Using knowledge of the true distributions in \eqref{eq: ex_MDP}, we compute the true parameters, $\theta_t^*$ for $t\in \{1,\cdots,T\}$, via the results in Appendix~\ref{subsec: true weights}, which we then use to compute the true optimal policy $\pi_t^*$ as in Corollary \ref{cor: opt_policy}. Finally, we use the true parameters, we compute the regret in \eqref{eq: regret}. We present our numerical results in Fig.~\ref{fig: example}. The regret $\mathcal{R}(L)$ as a function of the number of episodes $L$ is shown in Fig. \ref{fig: example}(a). We observe that the regret scales as $\widetilde{\mathcal{O}}(\sqrt{L})$, which is consistent with our results in Theorem~\ref{thrm: regret bound}. The average regret, $\mathcal{R}(L)/L$ as a function of $L$ is shown in Fig.~\ref{fig: example}(b) and is observed to converge as $L$ increases, indicating convergence of our algorithm. Fig~\ref{fig: example}(c) presents the norm of the estimation error between the learned and the true parameters, averaged over the episode horizon $T$, as a function of $L$. This is expressed as $\|\theta^{L}_{1:T}-\theta^{*}_{1:T}\|/T$, where $\theta^{L}_{1:T}\in \mathbb{R}^{d(n+1)\times T}$ is a matrix whose columns corresponds to the parameters $\theta^{L}_t$ at each time step $t$, and similarly for $\theta^{*}_{1:T}$. We observe that the estimation error decreases with $L$, indicating that the learned policy gradually converges to the optimal one. Finally, we apply both the learned policy at episode $L=1000$ and the optimal policy to the system in \eqref{eq: ex_dynamics}, and compare their corresponding closed-loop state trajectories, as shown in Fig.~\ref{fig: example}(d), alongside the trajectory of the exogenous state $s$. We observe that the trajectory under the learned policy closely matches that of the optimal policy, and both effectively track the mean of the exogenous state.

\section{Conclusion}\label{sec: conclusion}
In this work, we proposed a reinforcement learning framework that unifies linear control systems and feature-based linear Markov models, capturing both deterministic system dynamics and stochastic environmental effects. By leveraging this structure, we derived closed-form expressions forthe optimal value function and policy, and introduced a least-squares value iteration algorithm that learns the optimal control policy without requiring explicit model identification or exploration. We provided theoretical guarantees on stability and convergence, and demonstrated the effectiveness~of our approach through numerical simulations. Future directions include extending our method to settings with unknown system dynamics (with additive process noise), richer feature models, and more general classes of stochastic processes.

\bibliographystyle{unsrt}
\bibliography{alias,Main_bib}

\begin{thebibliography}{10}

\bibitem{BDA-JBM:07}
B.~D. Anderson and J.~B. Moore.
\newblock {\em Optimal control: linear quadratic methods}.
\newblock Courier Corporation, 2007.

\bibitem{SJB-AGB:96}
S.~J. Bradtke and A.~G. Barto.
\newblock Linear least-squares algorithms for temporal difference learning.
\newblock {\em Machine learning}, 22(1):33--57, 1996.

\bibitem{FSF-MIR:07}
F.~S. Francisco and M.~I. Ribeiro.
\newblock Q-learning with linear function approximation.
\newblock In {\em International Conference on Computational Learning Theory},
  pages 308--322. Springer, 2007.

\bibitem{RSS-AGB:18}
R.~S. Sutton and A.~G. Barto.
\newblock {\em Reinforcement learning: An introduction}.
\newblock MIT press, 2018.

\bibitem{CDP-PT:19}
C.~{De Persis} and P.~Tesi.
\newblock Formulas for data-driven control: Stabilization, optimality and
  robustness.
\newblock {\em IEEE Transactions on Automatic Control}, 65(3):909--924, 2020.

\bibitem{FD-PT-CDP:22}
F.~D{\"o}rfler, P.~Tesi, and C.~{De Persis}.
\newblock On the role of regularization in direct data-driven {LQR} control.
\newblock In {\em {IEEE} Conf.\ on Decision and Control}, pages 1091--1098,
  Canc\'un, Mexico, December 2022. IEEE.

\bibitem{FC-GB-FP:23a}
F.~Celi, G.~Baggio, and F.~Pasqualetti.
\newblock Closed-form and robust formulas for data-driven {LQ} control.
\newblock {\em {Annual Reviews in Control}}, 56, 2023.

\bibitem{WA-DK-BJ-RM-MS:05}
W.~Aangenent, D.~Kostic, B.~de~Jager, R.~van~de Molengraft, and M.~Steinbuch.
\newblock Data-based optimal control.
\newblock In {\em {A}merican {C}ontrol {C}onference}, pages 1460--1465,
  Portland, OR, USA, June 2005.

\bibitem{GRGDS-ASB-CL-LC:18}
G.~R.~G. da~Silva, A.~S. Bazanella, C.~Lorenzini, and L.~Campestrini.
\newblock Data-driven {LQR} control design.
\newblock {\em IEEE Control Systems Letters}, 3(1):180--185, 2018.

\bibitem{SD-HM-NM-BR-ST:20}
S.~Dean, H.~Mania, N.~Matni, B.~Recht, and S.~Tu.
\newblock On the sample complexity of the linear quadratic regulator.
\newblock {\em Foundations of Computational Mathematics}, 20(4):633--679, 2020.

\bibitem{MF-RG-SK-MM:18}
M.~Fazel, R.~Ge, S.~Kakade, and M.~Mesbahi.
\newblock Global convergence of policy gradient methods for the linear
  quadratic regulator.
\newblock In {\em International Conference on Machine Learning}, pages
  1467--1476, Stockholm, Sweden, 2018.

\bibitem{HM-AZ-MS-MRJ:19}
H.~Mohammadi, A.~Zare, M.~Soltanolkotabi, and M.~R. Jovanovi{\'c}.
\newblock Global exponential convergence of gradient methods over the nonconvex
  landscape of the linear quadratic regulator.
\newblock In {\em {IEEE} Conf.\ on Decision and Control}, pages 7474--7479,
  Nice, France, Dec. 2019.

\bibitem{JB-AM-MF-MM:19}
J.~Bu, A.~Mesbahi, M.~Fazel, and M.~Mesbahi.
\newblock Lqr through the lens of first order methods: Discrete-time case.
\newblock {\em arXiv preprint arXiv:1907.08921}, 2019.

\bibitem{IF-BP:21}
I.~Fatkhullin and B.~Polyak.
\newblock Optimizing static linear feedback: Gradient method.
\newblock {\em SIAM Journal on Control and Optimization}, 59(5):3887--3911,
  2021.

\bibitem{SJB-BEY-AGB:94}
S.~J. Bradtke, B.~E. Ydstie, and A.~G. Barto.
\newblock Adaptive linear quadratic control using policy iteration.
\newblock In {\em {A}merican {C}ontrol {C}onference}, volume~3, pages
  3475--3479. IEEE, 1994.

\bibitem{JK-JAB-JP:13}
J.~Kober, J.~A. Bagnell, and J.~Peters.
\newblock Reinforcement learning in robotics: A survey.
\newblock {\em The International Journal of Robotics Research},
  32(11):1238--1274, 2013.

\bibitem{VM-KK-DS-AG:13}
V.~Mnih, K.~Kavukcuoglu, D.~Silver, A.~Graves, I.~Antonoglou, D.~Wierstra, and
  M.~Riedmiller.
\newblock Playing atari with deep reinforcement learning.
\newblock {\em arXiv preprint arXiv:1312.5602}, 2013.

\bibitem{SD-HA-CJ:2016}
D.~Silver, A.~Huang, C.~J. Maddison, A.~Guez, L.~Sifre, G.~Van~Den Driessche,
  J.~Schrittwieser, I.~Antonoglou, V.~Panneershelvam, M.~Lanctot, and et~al.
\newblock Mastering the game of {G}o with deep neural networks and tree search.
\newblock {\em Nature}, 529(7587):484--489, 2016.

\bibitem{LY-MW:19}
L.~Yang and M.~Wang.
\newblock Sample-optimal parametric q-learning using linearly additive
  features.
\newblock In {\em Proceedings of the 36th International Conference on Machine
  Learning}, pages 6995--7004, 2019.

\bibitem{LY-MW:20}
L.~Yang and M.~Wang.
\newblock Reinforcement learning in feature space: Matrix bandit, kernels, and
  regret bound.
\newblock In {\em International Conference on Machine Learning}, pages
  10746--10756. PMLR, 2020.

\bibitem{CJ-ZY-ZW-MIJ:20}
C.~Jin, Z.~Yang, Z.~Wang, and M.~I. Jordan.
\newblock Provably efficient reinforcement learning with linear function
  approximation.
\newblock In Jacob Abernethy and Shivani Agarwal, editors, {\em Proceedings of
  Thirty Third Conference on Learning Theory}, volume 125 of {\em Proceedings
  of Machine Learning Research}, pages 2137--2143. PMLR, 09--12 Jul 2020.

\bibitem{NJ-AK-AA-JL-RES:17}
N.~Jiang, A.~Krishnamurthy, A.~Agarwal, J.~Langford, and R.~E. Schapire.
\newblock Contextual decision processes with low {B}ellman rank are
  {PAC}-learnable.
\newblock In {\em Proceedings of the 34th International Conference on Machine
  Learning}, pages 1704--1713, 2017.

\bibitem{RM:05}
R.~Munos.
\newblock Error bounds for approximate value iteration.
\newblock In {\em Proceedings of the National Conference on Artificial
  Intelligence}, volume~20, page 1006, 2005.

\bibitem{AZ-AL-MK-EB:20}
A.~Zanette, A.~Lazaric, M.~Kochenderfer, and E.~Brunskill.
\newblock Learning near optimal policies with low inherent bellman error.
\newblock In {\em International Conference on Machine Learning}, pages
  10978--10989. PMLR, 2020.

\bibitem{WS-NJ-AK-AA-JL:19}
W.~Sun, N.~Jiang, A.~Krishnamurthy, A.~Agarwal, and J.~Langford.
\newblock Model-based {RL} in contextual decision processes: Pac bounds and
  exponential improvements over model-free approaches.
\newblock In {\em Conference on learning theory}, pages 2898--2933. PMLR, 2019.

\bibitem{SD-SK-JL-Sl-GM-WS-RW:21}
S.~Du, S.~Kakade, J.~Lee, S.~Lovett, G.~Mahajan, W.~Sun, and R.~Wang.
\newblock Bilinear classes: A structural framework for provable generalization
  in rl.
\newblock In {\em International Conference on Machine Learning}, pages
  2826--2836. PMLR, 2021.

\bibitem{KO:10}
K.~Ogata.
\newblock {\em Modern Control Engineering}.
\newblock Instrumentation and controls series. Prentice Hall, 2010.

\bibitem{FC-GB-FP:22}
F.~Celi, G.~Baggio, and F.~Pasqualetti.
\newblock Closed-form estimates of the {LQR} gain from finite data.
\newblock In {\em {IEEE} Conf.\ on Decision and Control}, pages 4016--4021,
  Canc\'un, Mexico, December 2022.

\bibitem{JAT:12}
J.~A. Tropp.
\newblock User-friendly tail bounds for sums of random matrices.
\newblock {\em Foundations of Computational Mathematics}, 12(4):389--434, 2012.

\bibitem{YAY-DP-CS:11}
Y.~Abbasi-Yadkori, D.~P{\'a}l, and C.~Szepesv{\'a}ri.
\newblock Improved algorithms for linear stochastic bandits.
\newblock In {\em Advances in Neural Information Processing Systems}, 2011.

\end{thebibliography}

\appendices

\makeatletter
\renewcommand{\thesubsection}{\arabic{subsection}}
\renewcommand{\thesubsectiondis}{\thesection.\thesubsection.}
\setcounter{subsection}{0}
\renewcommand{\p@subsection}{\thesection.}
\makeatother

\section{Proof of Theorem \ref{thrm: Qfunction}}\label{app: A}




We first present the following result in which we provide an expression for the optimal greedy policy, $\pi_t^*$, and the optimal value function, $V_t^*(x,s)$ under the greedy policy.
\begin{theorem}{\bf \emph{(optimal policy and value function)}}\label{thrm: Vfunction}
Consider the dynamics in \eqref{eq: system} and the Markov Process in \eqref{eq: environment}. Let Assumption \ref{assump: linMDP} be satisfied. Then, for any $x\!\in\! \mc{X}$, $s\!\in\! \mc{S}$, $s'\!\sim\! \mathbb{P}(.|s)$, and~$t\!\geq0\!$
%
\begin{align}
 u^*(x,s,t)=\underbrace{K_{t,x}x + K_{t,h}\expect{h_{t+1}(s')|s} + K_{t,s}s   }_{\pi_t^*(x,s)},
 \end{align}
 where
 \begin{align}\label{eq: feedback gains in theorem}
\begin{split}
 K_{t,x}&=-\left(R +B^{\T} G_{t+1} B \right)^{-1} \left( B^{\T}G_{t+1}A+D^{\T}\right),\\
  K_{t,h}&=-\left(R +B^{\T} G_{t+1} B \right)^{-1}B^{\T},\\
 K_{t,s}&=-\left(R +B^{\T} G_{t+1} B \right)^{-1}H^{\T}.
 \end{split}
\end{align}
Further,
\begin{align}
 V_t^*(x,s)=x^{\T}G_t x + 2h_t^{\T}(s)x +q_t(s),
\end{align}
where $G_t \in \mathbb{R}^{n\times n} \succ 0$, $h_t(\cdot)\in \mathbb{R}^n$, and $q_t(\cdot)\in \mathbb{R}$ satisfy
\begin{align}
\begin{split} \label{eq: riccati_eq}
& G_t=A^{\T}G_{t+1}A + W \\
& \quad-\! (A^{\T}\!G_{t+1}B \!+\! D ) (R \!+\! B^{\T}G_{t+1}B )^{-1}\!(B^{\T}\!G_{t+1}A \!+\! D^{\T}),
 \end{split}
 \\
\begin{split}\label{eq: h}
&h_t\left(s_t\right)=\left(A^{\T} + K_{t,x}^{\T}B^{\T}\right)\Expect{}{h_{t+1}\left(s_{t+1}\right)}{s_t}\\
&\qquad \quad \quad+\left(F +K_{t,x}^{\T}H^{\T}\right)s_t, 
 \end{split}
 \\
\begin{split}\label{eq: q}
& q_{t}\left(s_t\right)=\Expect{}{q_{t+1}\left(s_{t+1}\right)}{s_t}+{s_t}^{\T}(M + H K_{s,t})s_t \\
 &\qquad \qquad+ \Expect{}{h_{t+1}^{\T}(s_{t+1})}{s_t} B K_{t,h} \Expect{}{h_{t+1}(s_{t+1})}{s_t}\\
 &\qquad \qquad+2 {s_t}^{\T}H K_{t,h}\Expect{}{h_{t+1}\left(s_{t+1}\right)}{s_t}, 
 \end{split}
\end{align}
with $G_T\!=\!W$, and $h_T(s_{T})=Fs_{T}$ and $q_{T}(s_{T})=s_{T}^{\T} M s_{T}$.
\begin{proof}
 We prove our claim by induction. For notational convenience, we drop the time index from the states and inputs inside the expressions and arguments of $c(\cdot)$, $V_t^*(\cdot)$, and $Q_t^*(\cdot)$ for $t\geq0$, where we use $x$, $s$, $u$, $x'$, and $s'$ to denote $x_t$, $s_t$, $u_t$, $x_{t+1}$, and $s_{t+1}$, respectively. At $t=T-1$,
\begin{align}\label{eq: Q(T-1)}
\begin{split}
&Q_{T-1}^*(x, s, u)\\\
&\quad=c(x,s,u) + \Expect{s'\sim \mathbb{P}_t(s'|s)}{V^*_{T}(x',s')}{x,s,u}\\
&\quad=x^{\T}(W\!+\!A^{\T}W A)x\! +\!u^{\T}\!(R \!+\!B^{\T} W  B)u \\
&\quad\quad+2(s^{\T}F^{\T} \!+\!{\Expect{}{s'}{s}}^{\T} F^{\T} A)x \\
&\quad \quad+2 \left(x^{\T}\left( D + A^{\T}W B\right) + s^{\T}H  + {\Expect{}{s'}{s}}^{\T}F^{\T} B\right)u\\
&\quad \quad+ s^{\T}M s +\Expect{}{s'^{\T}M s'}{s},
 \end{split}
\end{align}
where we used the fact that $V^*_T(x_T,s_T)=c(x_T,s_T,u_T)$ and $u_{T}=0$. Taking the derivative of $Q_{T-1}^*$ with respect~to~$u$
\begin{align*}
\frac{\partial Q_{T-1}^*(x,s,u)}{\partial u}=& 2(R +B^{\T}W B)u + 2(B^{\T}W A+D^{\T})x\\
& +2 (B^{\T} F\Expect{}{s'}{s} +  H^{\T}s). 
\end{align*}
Setting the above derivative to zero and solving for $u$, we get
\begin{align}\label{eq: u_{T-1}}
\begin{split}
 u^*_{T-1}&=-(R + B^{\T}W B  )^{-1} (B^{\T}W A + D^{\T})x_{T-1}\\
  &-(R + B^{\T}W B  )^{-1}(B^{\T} F \Expect{}{s_{T}}{s_{T-1}} +  H^{\T} s_{T-1}),
  \end{split}
\end{align} 
which is the minimizer of \eqref{eq: Q(T-1)}. We substitute \eqref{eq: u_{T-1}} in \eqref{eq: Q(T-1)},
\begin{align}\label{eq: V(T-1)}
V_{T-1}^*(x,s)\!=\!x^{\T}G_{T-1}x \!+\! 2 {h_{T-1}\left(s\right)}^{\T}x \!+\! q_{T-1}\left(s\right),
\end{align}
where $G_{T-1}$, $h_{T-1}\left(s\right)$, and $q_{T-1}\left(s\right)$ are as in Theorem \ref{thrm: Vfunction} for $t=T-1$. Suppose for $t=k+1$,
\begin{align*}
V_{k+1}^*(x,s)=x^{\T}G_{k+1}x + 2 {h_{k+1}\left(s\right)}^{\T}x + q_{k+1}\left(s\right), 
\end{align*}
where $G_{k+1}$, $h_{k+1}\left(s\right)$, and $q_{k+1}\left(s\right)$ are as in Theorem \ref{thrm: Vfunction} for $t=k+1$. Then we have,
\begin{align}\label{eq: Q(k)}
\begin{split}
&Q_{k}^*(x, s, u)\\
&\!\!\!\!\!\quad=c(x,s,u) + \Expect{s'\sim \mathbb{P}_t(s'|s)}{V^*_{k+1}(x',s')}{x,s,u}\\
&\!\!\!\!\!\quad=x^{\T}\!(W \!+\!A^{\T}G_{k+1} A)x\! +\!u^{\T}\!(R \!+\!B^{\T} G_{k+1} B)u \\
&\!\!\!\!\!\qquad+2(s^{\T}F^{\T} \!+\!{\Expect{}{h_{k+1}(s')}{s}}^{\T}A)x \\
&\!\!\!\!\!\qquad+\!2 \Big(x^{\T}( D \!+\! A^{\T}G_{k+1} B) \!+\! s^{\T}H \!+\! {\Expect{}{h_{k+1}(s')}{s}}^{\T}B\Big)u\\
&\!\!\!\!\!\qquad+ s^{\T}M s +\Expect{}{q_{k+1}(s')}{s},
\end{split}
\end{align}
Taking the derivative of $Q_{k}^*$ with respect~to $u$
\begin{align*}
\frac{\partial Q_{k}^*(x,s,u)}{\partial u}\!=& 2(R \!+ \!B^{\T}G_{k+1}B)u\!+\! 2(B^{\T}G_{k+1}A\!+\!D^{\T} )x \\
&+2\left(B^{\T} \Expect{}{h_{k+1}(s')}{s} +  H^{\T}s\right). 
\end{align*}
Setting the above derivative to zero and solving for $u$, we get
\begin{align}\label{eq: u(k)}
\begin{split}
 u^*_k=&-(R + B^{\T}G_{k+1}B  )^{-1} (B^{\T}G_{k+1}A+D^{\T})x_k\\
  &-(R + B^{\T}G_{k+1}B  )^{-1}B^{\T} \Expect{}{h_{k+1}(s(k+1))}{s_{k}} \\
  &-(R + B^{\T}G_{k+1}B  )^{-1}H^{\T}s_{k},
  \end{split}
\end{align} 
which is the minimizer of \eqref{eq: Q(k)}. We substitute \eqref{eq: u(k)} in \eqref{eq: Q(k)},
\begin{align*}
V_{k}^*(x,s)=x^{\T}G_{k}x + 2 {h_{k}\left(s\right)}^{\T}x + q_{k}\left(s\right),
\end{align*}
where $G_{k}$, $h_{k}\left(s\right)$, and $q_{k}\left(s\right)$ are as in Theorem \ref{thrm: Vfunction} for $t=k$. This completes the proof.
\end{proof}
\end{theorem}
\begin{pfof}{Theorem \ref{thrm: Qfunction}}
For notational convenience, we drop the time index from the states and inputs inside the expressions and arguments of $c(\cdot)$, $V_t^*(\cdot)$, and $Q_t^*(\cdot)$ for $t\geq0$, where we use $x$, $s$, $u$, $x'$, and $s'$ to denote $x_t$, $s_t$, $u_t$, $x_{t+1}$, and $s_{t+1}$, respectively. Using Theorem \ref{thrm: Vfunction}, we write
\begin{align*}
 &Q_t^*(x,s,u)\\
 &\quad =c(x,s,u) + \Expect{s'\sim \mathbb{P}_t(s'|s)}{V^*_{t+1}(x',s')}{x,s,u}\\
 &\quad=c(x,s,u) + \int_{\mc{S}}V^*_{t+1}(Ax+Bu,s') \mathbb{P}_t( \text{d}s'|s)\\
&\quad \overset{\text{(a)}}{=}c(x,s,u)+ \int_{\mc{S}}V^*_{t+1}(Ax+Bu,s') {\phi(s)}^{\T} \mu_t(\text{d}s')\\
&\quad \overset{\text{(b)}}{=}c(x,s,u) \! +\! \left(Ax\!+\!Bu\right)^{\T}\! G_{t+\!1}\! \left(Ax\!+\!Bu\right)\! \\
 &\quad \quad+ \int_{\mc{S}}(2 {h_{t+1}(s')}^{\T}(Ax\!+\!Bu) \!+\!q_{t+1}(s'))\\
 &\quad \qquad \qquad .\sum_{i=1}^d {\phi_i(s)} \mu_{i,t}(\text{d}s')\\
 &\quad=c(x,s,u) \! +\! \left(Ax\!+\!Bu\right)^{\T}\! G_{t+\!1}\! \left(Ax\!+\!Bu\right)\\ 
 &\quad \quad+\sum_{i=1}^d {\phi_i(s)} \int_{\mc{S}}q_{t+1}(s')\mu_{i,t}(\text{d}s')\\
 &\quad \quad+2\sum_{i=1}^d \Big( \phi_i(s)  \int_{\mc{S}}{h_{t+1}(s')}^{\T}\mu_{i,t}(\text{d}s')\Big)(Ax+Bu)\\
  &\quad =c(x,s,u) \! +\! \left(Ax\!+\!Bu\right)^{\T}\! G_{t+\!1}\! \left(Ax\!+\!Bu\right)\\
  &\quad \quad +\sum_{i=1}^d {\phi_i(s)} \underbrace{\Expectover{\mu_{i,t}}{q_{t+1}(s')}}_{\overline{q}_{i,t+1}}\\
 &\quad \quad+2\sum_{i=1}^d \Big( \phi_i(s)  \underbrace{\Expectover{\mu_{i,t}}{h_{t+1}(s')}^{\T}}_{\overline{h}^{\T}_{i,t+1}} \Big)(Ax+Bu),
\end{align*}
where in step (a) we have used Assumption \ref{assump: linMDP}, and in step (b) we have used Theorem \ref{thrm: Vfunction}.
\end{pfof}
%
%
%
%
%
%
%
%
\section{Least-squares value iteration}\label{app: LSVI}
We formulate the regularized least squares regression and derive its solution presented in lines $6$-$7$ of Algorithm \ref{alg: alg1}.  We begin by using the notation in \eqref{eq: theta} to derive the expression of the parametrized state-action value function, $Q_t$ and the corresponding parametrized optimal greedy policy in \eqref{eq: policy alg}. Using the expression of $Q_t$ in Theorem \ref{thrm: Qfunction}, we can write
\begin{align} \label{eq: Qfunc proof} 
\begin{split}
 Q_t(x,s,u)=& c(x,s,u)\!+\! \left(Ax+Bu\right)^{\T} G_{t+1}\left(Ax+Bu\right)\\
 &+\sum_{i=1}^{d}\phi_i(s)\Big(2\left(Ax\!+\!Bu\right)^{\T} \overline{h}_{i,t+1}+\overline{q}_{i,t+1}\Big)\\
 =&c(x,s,u)\!+\! \left(Ax+Bu\right)^{\T} G_{t+1}\left(Ax+Bu\right)\\
 &+\sum_{i=1}^{d}\phi_i(s)\underbrace{\begin{bmatrix}
2 \left(Ax + Bu\right)^{\T} & 1
\end{bmatrix}}_{{y(x,u)}^{\T}}
\underbrace{\begin{bmatrix}
 \overline{h}_{i,t+1}\\
 \overline{q}_{i,t+1}
\end{bmatrix}}_{\theta_{i,t+1}}\\
=& c(x,s,u)\!+\! \left(Ax+Bu\right)^{\T} G_{t+1}\left(Ax+Bu\right)\\
 &+\sum_{i=1}^{d} \phi_{i}(s){y(x,u)}^{\T} \theta_{i,t+1} \\
 =& c(x,s,u)\!+\! \left(Ax+Bu\right)^{\T} G_{t+1}\left(Ax+Bu\right)\\
&+  {\phi(s)}^{\T} \underbrace{\left( I_d \otimes {y(x,u)}^{\T} \right)}_{Y(x,u)}
\underbrace{\begin{bmatrix}
 \theta_{1,t+1}\\
 \vdots \\
 \theta_{d,t+1}
\end{bmatrix}}_{\theta_{t+1}}\\
 =& c(x,s,u)\!+\! \left(Ax+Bu\right)^{\T} G_{t+1}\left(Ax+Bu\right)\\
&+ {\phi(s)}^{\T} Y(x,u) {\theta}_{t+1}.
 \end{split}
\end{align}
Next, using the notation in \eqref{eq: theta}, we can write
\begin{align}\label{eq: aux3}
\begin{split}
 \overline{h}_{i,t+1}=&
\underbrace{\begin{bmatrix}
 I_n & 0_{n\times 1}
\end{bmatrix}}_{Z}
\underbrace{\begin{bmatrix}
 \overline{h}_{i,t+1}\\
 \overline{q}_{i,t+1}
\end{bmatrix}}_{\theta_{i,t+1}}
=Z\theta_{i,t+1}, ~~ i\!\in\!\{1,\cdots,d\},\\
\overline{h}_{t+1}=&
\begin{bmatrix}
 Z & 0 & \cdots & 0\\
 0 & Z & \cdots & 0\\
 \vdots & \ddots & \ddots & \vdots\\
 0 & 0 & \cdots & Z
\end{bmatrix}
\begin{bmatrix}
 \theta_{1,t+1}\\
 \vdots\\
 \theta_{d,t+1}
\end{bmatrix}
=\left(I_d \otimes Z\right) \theta_{t+1}.
\end{split}
\end{align}
Similarly, we can write
\begin{align}\label{eq: aux7}
\begin{split}
\overline{q}_{i,t+1}=&
\underbrace{\begin{bmatrix}
 0_{1\times n} & 1
\end{bmatrix}}_{\overline{Z}}
\begin{bmatrix}
 \overline{h}_{i,t+1}\\
 \overline{q}_{i,t+1}
\end{bmatrix}
=\overline{Z}\theta_{i,t+1},, ~~ i\!\in\!\{1,\cdots,d\},\\\\
\overline{q}_{t+1}=&
\begin{bmatrix}
 \overline{Z} & 0 & \cdots & 0\\
 0 & \overline{Z} & \cdots & 0\\
 \vdots & \ddots & \ddots & \vdots\\
 0 & 0 & \cdots & \overline{Z}
\end{bmatrix}
\begin{bmatrix}
 \theta_{1,t+1}\\
 \vdots\\
 \theta_{d,t+1}
\end{bmatrix}
=\left(I_d \otimes \overline{Z}\right) \theta_{t+1}.
\end{split}
\end{align}
Next, from the expression of $u_t$ in Corollary \ref{cor: opt_policy}, we write
\begin{align}\label{eq: greedy policy proof}
\begin{split}
u_t(x,s)=&K_{t,x}x + K_{t,s}s + K_{t,h}\sum_{i=1}^{d} \phi_i(s) \overline{h}_{i,t+1}\\
=& K_{t,x}x + K_{t,s}s\\
&+ K_{t,h}
\begin{bmatrix}
 \phi_1(s) I_n & \cdots & \phi_d(s) I_n
\end{bmatrix}
\begin{bmatrix}
 \overline{h}_{1,t+1}\\
 \vdots\\
 \overline{h}_{d,t+1}
\end{bmatrix}\\
=& K_{t,x}x + K_{t,s}s+ K_{t,h}\left({\phi(s)}^{\T} \otimes I_n\right) \overline{h}_{t+1}\\
\overset{\text{(a)}}{=}&K_{t,x}x + K_{t,s}s+\left({\phi(s)}^{\T} \otimes I_n\right) \left(I_d \otimes Z\right) \theta_{t+1}\\
=&K_{t,x}x + K_{t,s}s+\left({\phi(s)}^{\T} \otimes Z\right) \theta_{t+1},
\end{split}
\end{align}
where in step (a) we have used \eqref{eq: aux3}. We define the Bellman target at time $t$ as
\begin{align}\label{eq: target1}
 g_t(x,s,u)=c(x,s,u) + \min_{v} \widehat{Q}_{t+1}(x',s',v),
\end{align}
where $x'$ and $s'$ denote the states resulting from taking action $u$ in states $x$ and $s$, and $\widehat{Q}_{t+1}(x',s',v)$ is the estimate of the state-action value function at time $t+1$. We re-write \eqref{eq: target1} as
\begin{align}\label{eq: target1}
\begin{split}
 g_t(x,s,u)=&c(x,s,u) + \widehat{V}^*_{t+1}(x',s')\\
 \overset{\text{(b)}}{=}&c(x,s,u)+ x'^{\T}G_{t+1} x'\\
  &+ 2\widehat{h}_{t+1}^{\T}(s')x' +\widehat{q}_{t+1}(s'),
  \end{split}
\end{align}
where in step (b), we have used Theorem \ref{thrm: Vfunction}. For notational convenience, let $X_1(t)\!=\!A^{\T} \!+\! K_{t,x}^{\T}B^{\T}$, $X_2(t)\!=\!F \!+\!K_{t,x}^{\T}H^{\T}$, $Y_1(t)\!=\!M \!+\! H K_{s,t}$, $Y_2(t)\!=\!B K_{t,h}$, and $Y_3(t)\!=\!H K_{t,h}$. Using \eqref{eq: h} and \eqref{eq: q}, we re-write $\widehat{h}_{t+1}$ and $\widehat{q}_{t+1}$ in \eqref{eq: target1}~as
\begin{align}\label{eq: h in terms of theta}
\begin{split}
 \widehat{h}_{t+1}\left(s_{t+1}\right)=&X_1(t+1)\Expect{}{h_{t+2}\left(s_{t+2}\right)}{s_{t+1}}\\
 &+X_2(t+1)s_{t+1},\\
 \overset{\text{(c)}}{=}&X_1(t+1)\left({\phi(s_{t+1})}^{\T} \otimes Z\right) \widehat{\theta}_{t+2} \\
 &+ X_2(t+1)s_{t+1},
 \end{split}
 \end{align}
 \begin{align}\label{eq: q in terms of theta}
\begin{split}
 & \widehat{q}_{t+1}\left(s_{t+1}\right)\\
 &~ =\Expect{}{q_{t+2}\left(s_{t+2}\right)}{s_{t+1}}+{s_{t+1}}^{\T}Y_1(t+1)s_{t+1}\\
 & ~+ \Expect{}{h_{t+2}^{\T}(s_{t+2})}{s_{t+1}} Y_2(t\!+\!1) \Expect{}{h_{t+2}(s_{t+2})}{s_{t+1}}\\
 & ~+2 {s_{t+1}}^{\T}Y_3(t+1)\Expect{}{h_{t+2}(s_{t+2})}{s_{t+1}},\\
\overset{\text{(d)}}{=}&{\phi(s_{t+1})}^{\T}\left(I_d \otimes \overline{Z}\right) \widehat{\theta}_{t+2} +{s_{t+1}}^{\T}Y_1(t+1)s_{t+1}\\
 &+ {\widehat{\theta}_{t+2}}^{\T} \left({\phi(s_{t+1})}\otimes Z^{\T}\right) Y_2(t\!+\!1) \left({\phi(s_{t+1})}^{\T} \otimes Z\right) \widehat{\theta}_{t+2}\\
 &+2 {s_{t+1}}^{\T}Y_3(t+1)\left({\phi(s_{t+1})}^{\T} \otimes Z\right) \widehat{\theta}_{t+2},
 \end{split}
\end{align}
where in steps (c) and (d) we have used \eqref{eq: aux3} and \eqref{eq: aux7}, respectively. The temporal difference (TD) error is written as
\begin{align}\label{eq: TD error}
\begin{split}
\varepsilon_t(x,s,u) =&g_t(x,s,u) - \widehat{Q}_{t}(x,s,u) \\
=& c(x,s,u) +\min_{v} \widehat{Q}_{t+1}(x',s',v)- \widehat{Q}_{t}(x,s,u) \\
\overset{\text{(d)}}{=}&c(x,s,u)\!+\! x'^{\T}G_{t+1} x'\!+\! 2\widehat{h}_{t+1}^{\T}(s')x' \!+\!\widehat{q}_{t+1}(s')\\
  &- c(x,s,u)\!-\! \left(Ax+Bu\right)^{\T} G_{t+1}\left(Ax+Bu\right)\\
&- {\phi(s)}^{\T} Y(x,u) \widehat{{\theta}}_{t+1}\\
=&2\widehat{h}_{t+1}^{\T}(s')x' \!+\!\widehat{q}_{t+1}(s')- {\phi(s)}^{\T} Y(x,u) \widehat{{\theta}}_{t+1},
\end{split}
\end{align}
where in step (d) we have used \eqref{eq: Qfunc proof} and \eqref{eq: target1}. The TD error $\varepsilon_t(x,s,u)$ in \eqref{eq: TD error} captures the discrepancy between the Bellman target and the current estimate of the $Q$-function. In Least-Squares Value Iteration (LSVI) in Algorithm \ref{alg: alg1}, we minimize the squared TD error over the dataset, \eqref{eq: dataset}, collected up to episode $\ell-1$, to obtain an updated estimate of the $Q$-function. Specifically, the parameters $\widehat{\theta}_t$ at episode $\ell$ denoted by $\theta^{\ell}_t$ is obtained by solving the following regularized least-squares problem
\begin{align}\label{eq: least squares}
\begin{split}
\theta^{\ell}_{t+1}=& \arg \min_{\theta}~\underbrace{\sum_{j=1}^{\ell-1} {\varepsilon_t(x^{j},s^{j},u^{j})}^{2} + \lambda {{\|\theta\|}_2^2}}_{J}\\
=&\arg \min_{\theta} \sum_{j=1}^{\ell-1} \Big( {\phi(s^j)}^{\T} Y(x^j,u^j) \widehat{{\theta}}_{t+1} \\
&\qquad\qquad  ~~~ -2\widehat{h}_{t+1}^{\T}(s'^j)x'^j \!-\!\widehat{q}_{t+1}(s'^j)\Big)^2 \!+\! \lambda {{\|\theta\|}_2^2}.
 \end{split}
\end{align}
Taking the derivative of \eqref{eq: least squares} with respect to $\theta$, we get
\begin{align*}
 \frac{\partial J}{\partial \theta}\! =& 2\sum_{j=1}^{\ell-1} Y^{\T}(x^{j}\!,\!u^{j})\phi(s^{j})\Big( {\phi(s^{j})}^{\T} Y(x^{j},u^{j}) \theta\\
 &\!\!\!\qquad \qquad  -2\widehat{h}_{t+1}^{\T}(s'^j)x' \!-\!\widehat{q}_{t+1}(s'^j)\Big)+2\lambda \theta
\end{align*}
Setting the above derivative to zero and solving for $\theta$, we get
\begin{align}\label{eq: theta_star}
\begin{split}
 \theta_{t+1}^{\ell}\!=&\Lambda_t^{-1} \!\sum_{j=1}^{\ell-1}\!Y^{\T}\!(x^{j},u^{j})\phi(s^{j}) \big(2\widehat{h}_{t+1}^{\T}(s'^j)x'^j \!+\!\widehat{q}_{t+1}(s'^j) \big),\\
 \Lambda_t =& \sum_{j=1}^{\ell-1} {Y(x^{j},u^{j})}^{\T}\phi(s^{j}) {\phi(s^{j})}^{\T} Y(x^{j},u^{j}) + \lambda I_{d(n+1)}.
 \end{split}
\end{align}
In Algorithm \ref{alg: alg1} we computed $\widehat{h}_{t+1}$ and $\widehat{q}_{t+1}$ at episode $\ell$ as in \eqref{eq: h in terms of theta} and \eqref{eq: q in terms of theta}, respectively. These quantities are obtained using the updated parameter $\theta^{\ell}_{t+2}$ from the previous iteration of the backward-in-time weight update loop(lines 5-11 in Algorithm~\ref{alg: alg1}). In particular, we have $\widehat{h}_{t+1}(\cdot)\!=\!h^{\ell}_{t+1}(\cdot)$ and $\widehat{q}_{t+1}(\cdot)\!=\!~q^{\ell}_{t+1}(\cdot)$.

\section{True weights $\overline{h}_t$ and $\overline{q}_t$}\label{app: true theta}
Throughout this Appendix, we use the following notation, 
\begin{align}\label{eq: notations}
 \begin{split}
& X_1(t)\!=A^{\T} \!+\! K_{t,x}^{\T}B^{\T},\quad X_2(t)\!=\!F \!+\!K_{t,x}^{\T}H^{\T},\\
& Y_1(t)\!=M \!+\! H K_{s,t}, \quad Y_2(t)\!=\!B K_{t,h}, \quad Y_3(t)\!=\!H K_{t,h},\\
 &\Phi_t=
  \begin{bmatrix}
\Expectover{\mu_{1,t}}{\phi(s_t)^{\T}} \\
 \vdots \\
 \Expectover{\mu_{d,t}}{\phi(s_t)^{\T}}
\end{bmatrix},
~ \text{and} ~~
\overline{m}_t\!=\!
\begin{bmatrix}
\overline{m}_{1,t}\\
 \vdots \\
 \overline{m}_{d,t}
\end{bmatrix}
\!=\!\! \begin{bmatrix}
\Expectover{\mu_{1,t}}{s_t} \\
 \vdots \\
 \Expectover{\mu_{d,t}}{s_t}
\end{bmatrix}\!\!,
\end{split}
\end{align}
for $t\!\in\!\{0,\!\cdots\!,T\}$. In addition, we define for $t\!\in\!\{0,\!\cdots,\!T\}$
\begin{align}\label{eq: notations bounds}
\begin{split}
 \|X_1(t)\|\leq &\overline{X}_1, \quad  \|X_2(t)\|\leq \overline{X}_2,\\
 \quad\|Y_1(t)\|\leq &\overline{Y}_1, \quad\|Y_2(t)\|\leq \overline{Y}_2, \quad\|Y_3(t)\|\leq \overline{Y}_3.
\end{split}
\end{align}


\subsection{Closed-form expressions of $\overline{h}_t$ and $\overline{q}_t$}\label{subsec: true weights}
In this Appendix, we derive closed-form expressions for the true parameters $\overline{h}_t\!=\!\Expectover{\mu_{t}}{h_{t}(s_t)}$ and $\overline{q}_t\!=\!\Expectover{\mu_{t}}{q_{t}(s_t)}$.

\begin{theorem}{\bf \emph{(closed-form expressions for the true $\overline{h}_t$ and $\overline{q}_t$)}}\label{thrm: true parameters} Consider the dynamics in \eqref{eq: system} and the Markov Process in \eqref{eq: environment}. Let Assumption \ref{assump: linMDP} and Assumption \ref{assump: QuadCost} be satisfied. Let $\overline{h}_t\!=\!\Expectover{\mu_{t}}{h_{t}(s_t)}$ and $\overline{q}_t\!=\!\Expectover{\mu_{t}}{q_{t}(s_t)}$ for $t\in \{0,\cdots,T\}$, where $h_t(\cdot)$ and $q_t(\cdot)$ are as in \eqref{eq: h} and \eqref{eq: q}, respectively. Then, for $t \in\{0,\cdots, T-1\}$
\begin{align*}
 \begin{split}
 \overline{h}_t=& (\Phi_t \otimes X_1(t)) \overline{h}_{t+1} + (I_d \otimes X_2(t)) \overline{m}_t,\\
\overline{q}_t
=&\Phi_t \overline{q}_{t+1} +\Expectover{\mu_t}{s_t^{\T}Y_1(t) s_t}\\
&+ \begin{bmatrix}
 \overbar{h}_{t+1}^{\T}\! \left(\Expectover{\mu_{1,t}}{\phi(s_t)\phi(s_t)^{\T}} \!\otimes\! Y_2(t)\right)\!\overbar{h}_{t+1}\\
 \vdots \\
\overbar{h}_{t+1}^{\T} \!\left(\Expectover{\mu_{d,t}}{\phi(s_t)\phi(s_t)^{\T}}\! \otimes\! Y_2(t)\right)\!\overbar{h}_{t+1}
\end{bmatrix}\\
&+2 \Expectover{\mu_t}{\phi(s_t)^{\T} \otimes s_t^{\T}Y_3(t)}\overline{h}_{t+1},
\end{split}
\end{align*}
with $\overline{h}_T\!=\!F \Expectover{\mu_T}{s_T}$ and $\overline{q}_T=\Expectover{\mu_T}{s_T^{\T}M s_T}$, where $X_1(t)$, $X_2(t)$, $Y_1(t)$, $Y_2(t)$, $Y_3(t)$, $\Phi_t$, and $\overline{m}_t$ are as in \eqref{eq: notations}.
\begin{proof}
We re-write equation \eqref{eq: h} as
\begin{align}
h_t\left(s_t\right)=&X_1(t)\Expect{}{h_{t+1}\left(s_{t+1}\right)}{s_t}+X_2(t) s_t.
\end{align}
Taking the expectation of both sides with respect to $\mu_{i,t}$ for each $i\in\{1,\cdots , d\}$, we get
\begin{align}
\begin{split}
 &\Expectover{\mu_{i,t}}{h_t\left(s_t\right)}\\
 &\quad=X_1(t)\Expectover{\mu_{i,t}}{\Expect{}{h_{t+1}\left(s_{t+1}\right)}{s_t}}+X_2(t) \Expectover{\mu_{i,t}}{s_t}\\
 &\quad=X_1(t) \Expectover{\mu_{i,t}}{\sum_{j=1}^d  \phi_j(s_t) \Expectover{\mu_{j,t+1}}{h_{t+1}(s_{t+1})} }\\
  &\qquad+ X_2(t) \overline{m}_{i,t}\\
  &\quad=X_1(t) \sum_{j=1}^d\Expectover{\mu_{i,t}}{  \phi_j(s_t)} \Expectover{\mu_{j,t+1}}{h_{t+1}(s_{t+1})}  \\
  &\qquad+ X_2(t) \overline{m}_{i,t}\\
  %
%
&\quad= X_1(t) (\Expectover{\mu_{i,t}}{\phi(s_t)^{\T}} \otimes I_n) \overline{h}_{t+1} +  X_2(t) \overline{m}_{i,t}.
 \end{split}
\end{align}
By noting that
\begin{align*}
 \overline{h}_t=\begin{bmatrix}
 \Expectover{\mu_{1,t}}{h_{t}(s_t)}\\
 \vdots \\
 \Expectover{\mu_{d,t}}{h_{t}(s_t)}
\end{bmatrix},
\end{align*}
and denoting
\begin{align*}
 \Phi_t=
  \begin{bmatrix}
\Expectover{\mu_{1,t}}{\phi(s_t)^{\T}} \\
 \vdots \\
 \Expectover{\mu_{d,t}}{\phi(s_t)^{\T}}
\end{bmatrix},
\quad \text{and} \quad
\overline{m}_t=
\begin{bmatrix}
\overline{m}_{1,t}\\
 \vdots \\
 \overline{m}_{d,t}
\end{bmatrix},
\end{align*}
we can write
\begin{align}\label{eq: h_bar dynamics}
 \begin{split}
 \overline{h}_t=&(I_d \otimes X_1(t))(\Phi_t \otimes I_n) \overline{h}_{t+1} + (I_d \otimes X_2(t)) \overline{m}_t\\
 =& (\Phi_t \otimes X_1(t)) \overline{h}_{t+1} + (I_d \otimes X_2(t)) \overline{m}_t.
\end{split}
\end{align}
%
%
%
%
Next, we re-write equation \eqref{eq: q} as
\begin{align}
\begin{split}
 q_{t}\left(s_t\right)=&\Expect{}{q_{t+1}\left(s_{t+1}\right)}{s_t}+{s_t}^{\T}Y_1(t)s_t \\
 &+ \Expect{}{h_{t+1}^{\T}\left(s_{t+1}\right)}{s_t} Y_2(t) \Expect{}{h_{t+1}\left(s_{t+1}\right)}{s_t}\\
&+2 {s_t}^{\T}Y_3(t)\Expect{}{h_{t+1}\left(s_{t+1}\right)}{s_t}.
\end{split}
\end{align}
Taking the expectation of both sides with respect to $\mu_{i,t}$ for each $i\in\{1,\cdots , d\}$, we get
\begin{align}\label{eq: expect q_i}
\begin{split}
 &\Expectover{\mu_{i,t}}{q_{t}(s_t)}\!=\! \Expectover{\mu_{i,t}}{\Expect{}{q_{t+1}(s_{t+1})}{s_t}} \!+\!  \Expectover{\mu_{i,t}}{{s_t}^{\T}Y_1(t)s_t}\\ 
&\qquad \quad +  \Expectover{\mu_{i,t}}{\Expect{}{h_{t+1}^{\T}(s_{t+1})}{s_t} Y_2(t) \Expect{}{h_{t+1}(s_{t+1})}{s_t}}\\
 &\qquad \quad+2 \Expectover{\mu_{i,t}}{ {s_t}^{\T}Y_3(t)\Expect{}{h_{t+1}(s_{t+1})}{s_t}}.
\end{split}
\end{align}
We start with,
\begin{align}\label{eq: term1}
\begin{split}
&\Expectover{\mu_{i,t}}{\Expect{}{q_{t+1}\left(s_{t+1}\right)}{s_t}} \\
&~= \Expectover{\mu_{i,t}}{\sum_{j=1}^d \phi_j(s_t) \Expectover{\mu_{j,t+1}}{q_{t+1}(s_{t+1})}}\\
&~= \Expectover{\mu_{i,t}}{\sum_{j=1}^d \phi_j(s_t)} \Expectover{\mu_{j,t+1}}{q_{t+1}(s_{t+1})}\\
&~=\!\!
\underbrace{\begin{bmatrix}
 \Expectover{\mu_{i,t}\!\!}{  \phi_1(s_t)} &\!\! \!\cdots\!\! \!& \Expectover{\mu_{i,t}\!\!}{\phi_d(s_t)}
\end{bmatrix}}_{\Expectover{\mu_{i,t}}{\phi(s_t)^{\T}}}\!\!
\underbrace{\begin{bmatrix}
 \Expectover{\mu_{1,t+1}\!\!}{q_{t+1}(s_{t+1})}\\
 \vdots \\
 \Expectover{\mu_{d,t+1}\!\!}{q_{t+1}(s_{t+1})}
\end{bmatrix}}_{\overline{q}_{t+1}}\!\!.
\end{split}
\end{align}
Next we have,
\begin{align}\label{eq: term2}
\begin{split}
& \Expectover{\mu_{i,t}}{\Expect{}{h_{t+1}^{\T}\left(s_{t+1}\right)}{s_t} Y_2(t) \Expect{}{h_{t+1}\left(s_{t+1}\right)}{s_t}}\\
&\quad  =\mathbb{E}_{\mu_{i,t}}\Bigg[\sum_{j=1}^d \phi_j(s_t)\Expectover{\mu_{j,t+1}}{h_{t+1}(s_{t+1})^{\T} } Y_2(t)\\
&\qquad \qquad  \qquad \qquad .\sum_{k=1}^d \phi_k(s_t)\Expectover{\mu_{k,t+1}}{h_{t+1}(s_{t+1})}\Bigg]\\
& \qquad = \mathbb{E}_{\mu_{i,t}}\Biggl[
\begin{bmatrix}
 \Expectover{\mu_{1,t+1}}{h_{t+1}(s_{t+1})} \\
  \vdots \\
 \Expectover{\mu_{d,t+1}}{h_{t+1}(s_{t+1})}
\end{bmatrix}^{\T}
\underbrace{\begin{bmatrix}
 \phi_1(s_t) I_n \\
 \vdots \\
 \phi_d(s_t) I_n
\end{bmatrix}}_{\phi(s_t) \otimes I_n}
Y_2(t)\\
&\qquad \qquad  \qquad \quad .\begin{bmatrix}
 \phi_1(s_t) I_n \\
 \vdots \\
 \phi_d(s_t) I_n
\end{bmatrix}^{\T}
\underbrace{\begin{bmatrix}
 \Expectover{\mu_{1,t+1}}{h_{t+1}(s_{t+1})} \\
  \vdots \\
 \Expectover{\mu_{d,t+1}}{h_{t+1}(s_{t+1})}
\end{bmatrix}}_{\overline{h}_{t+1}}
\Bigg]\\
&\qquad = \overbar{h}_{t+1}^{\T} \Expectover{\mu_{i,t}}{\left(\phi(s_t) \otimes I_n \right) Y_2(t) \left(\phi(s_t)^{\T} \otimes I_n \right)}\overbar{h}_{t+1}\\
&\qquad =\! \overbar{h}_{t+\!1}^{\T} \Expectover{\mu_i}{\!(\phi(s_t)\! \otimes \!I_n )\! (1\!\otimes \!Y_2(t)) (\phi(s_t)^{\T} \!\otimes\! I_n )}\overbar{h}_{t+1}\\
&\qquad = \overbar{h}_{t+1}^{\T} \left(\Expectover{\mu_i}{\phi(s_t)\phi(s_t)^{\T}} \otimes Y_2(t)\right)\overbar{h}_{t+1}.
\end{split}
\end{align}
Finally, we have
\begin{align}\label{eq: term3}
\begin{split}
&\Expectover{\mu_{i,t}}{ {s_t}^{\T}Y_3(t)\Expect{}{h_{t+1}\left(s_{t+1}\right)}{s_t}} \\
&= \Expectover{\mu_{i,t}}{s_t^{\T} Y_3(t) \sum_{j=1}^d \phi_j (s_t)  \Expectover{\mu_{j,t+1}}{h_{t+1}(s_{t+1})} }\\
%
& =  \Expectover{\mu_{i,t}}{s_t^{\T} Y_3(t) \left(\phi(s_t)^{\T} \otimes I_n \right)} \overline{h}_{t+1}\\
& =  \Expectover{\mu_{i,t}}{ \left(1\otimes s_t^{\T}Y_3(t)\right) \left(\phi(s_t)^{\T} \otimes I_n \right)}\overline{h}_{t+1}\\
& =  \Expectover{\mu_{i,t}}{ \left(\phi(s_t)^{\T} \otimes s_t^{\T}Y_3(t) \right)}\overline{h}_{t+1}.
\end{split}
\end{align}
Substituting \eqref{eq: term1}, \eqref{eq: term2}, and \eqref{eq: term3} in \eqref{eq: expect q_i}, we get
\begin{align}
 \begin{split}
 \Expectover{\mu_{i,t}}{q_t(s_t)}=& \Expectover{\mu_{i,t}}{\phi(s_t)^{\T}} \overline{q}_{t+1} + \Expectover{\mu_{i,t}}{s_t^{\T}Y_1(t) s_t}\\
 &+ \overbar{h}_{t+1}^{\T} \left(\Expectover{\mu_{i,t}}{\phi(s_t)\phi(s_t)^{\T}} \otimes Y_2(t)\right)\overbar{h}_{t+1}\\
 &+2 \Expectover{\mu_{i,t}}{\phi(s_t)^{\T} \otimes s_t^{\T}Y_3(t)} \overline{h}_{t+1}.
\end{split}
\end{align}
Then, we can write
\begin{align}\label{eq: q_bar dynamics}
\begin{split}
 \underbrace{\begin{bmatrix}
 \Expectover{\mu_{1,t}}{q_{t}(s_t)}\\
 \vdots \\
 \Expectover{\mu_{d,t}}{q_{t}(s_t)}
\end{bmatrix}}_{\overline{q}_t}
&\!=\!
\underbrace{\begin{bmatrix}
 \Expectover{\mu_{1,t}\!\!}{\phi(s_t)^{\T}}\\
 \vdots \\
  \Expectover{\mu_{1,t}\!\!}{\phi(s_t)^{\T}}
\end{bmatrix}}_{\Phi_t}\!
 \underbrace{\begin{bmatrix}
 \Expectover{\mu_{1,t+1}\!\!}{q_{t+1}(s_{t+1})}\\
 \vdots \\
 \Expectover{\mu_{d,t+1}\!\!}{q_{t+1}(s_{t+1})}
\end{bmatrix}}_{\overline{q}_{t+1}}\\
&+\begin{bmatrix}
\Expectover{\mu_{1,t}}{s_t^{\T}Y_1(t) s_t}\\
 \vdots \\
\Expectover{\mu_{d,t}}{s_t^{\T}Y_1(t) s_t}
\end{bmatrix}\\
&\!\!\!+\! \begin{bmatrix}
 \overbar{h}_{t+1}^{\T}\! \left(\Expectover{\mu_{1,t}\!\!}{\phi(s_t)\phi(s_t)^{\T}} \!\otimes\! Y_2(t)\right)\!\overbar{h}_{t+1}\\
 \vdots \\
\overbar{h}_{t+1}^{\T} \!\left(\Expectover{\mu_{d,t}\!\!}{\phi(s_t)\phi(s_t)^{\T}}\! \otimes\! Y_2(t)\right)\!\overbar{h}_{t+1}
\end{bmatrix}\\
&+2 \begin{bmatrix}
 \Expectover{\mu_{1,t}\!\!}{\phi(s_t)^{\T} \otimes s_t^{\T}Y_3(t)}\\
 \vdots \\
\Expectover{\mu_{d,t}\!\!}{\phi(s_t)^{\T} \otimes s_t^{\T}Y_3(t)}
\end{bmatrix} \overline{h}_{t+1}.
\end{split}
\end{align}
Finally, from Theorem \ref{thrm: Vfunction}, we have $\overline{h}_T\!=\!\Expectover{\mu_{T}}{h_{T}(s_T)}=F\Expectover{\mu_{T}}{s_T}$, and $\overline{q}_t\!=\!\Expectover{\mu_{T}}{q_{T}(s_T)}=\Expectover{\mu_{T}}{s_T^{\T}M s_T}$.
\end{proof}
\end{theorem}

\subsection{Upper bounds on $\|\overline{h}_t\|$ and $\|\overline{q}_t\|$}\label{subsec: R_theta}
\begin{theorem}{\bf \emph{(upper bounds on the true $\overline{h}_t$ and $\overline{q}_t$)}}\label{thrm: true parameters bounds} Consider the dynamics in \eqref{eq: system} and the Markov Process in \eqref{eq: environment}. Let Assumption \ref{assump: linMDP} and Assumption \ref{assump: QuadCost} be satisfied. Let $\overline{h}_t$ and $\overline{q}_t$ be as in Theorem \ref{thrm: true parameters} for $t\in \{0,\cdots,T\}$. Then, 
\begin{align*}
 \begin{split}
  \|\overline{h}_t \| &\leq\left\|F\right\| \delta_s \alpha \rho^{T-t} +   \frac{\overline{X}_2 \delta_s \alpha \sqrt{d}}{1-\rho},\\
  \|\overline{q}_t\|&\leq \delta_s^2 \|M\| + \delta_s^2 \overline{Y}_1 \sqrt{d} + \frac{\overline{Y}_2 {\|\overline{h}_{t+1}\|}^2}{\sqrt{d}} + \delta_s \overline{Y}_3 \|\overline{h}_{t+1}\|,
\end{split}
\end{align*}
 for $t\in\{0,\cdots, T\}$ with $\| \overline{h}_{T+1}\|\!=\!0$, where $\alpha\!>\!0$, $0<\rho<1$ are constants, and $\overline{X}_1$, $\overline{X}_2$, $\overline{Y}_1$, $\overline{Y}_2$, and $\overline{Y}_3$ are as in~\eqref{eq: notations bounds}.
\begin{proof}
From \eqref{eq: h_bar dynamics}, we have
\begin{align}
\overline{h}_t =& (\Phi_t \otimes X_1(t)) \overline{h}_{t+1} + (I_d \otimes X_2(t)) \overline{m}_t.
\end{align}
Let $\Xi(t_1,t_2)=\prod_{i=t_2-1}^{t_1} \Phi_i \otimes X_1(i)$ with $t_2>t_1$. Then, given $\overline{h}_T$, we can write 
\begin{align}\label{eq: h_bar at t}
\begin{split}
 \overline{h}_t  &= \Xi(t,T)\overline{h}_T + \sum_{j=T-1}^{t} \Xi(t,j) \left(I_d \otimes X_2(j)\right) \overline{m}_j.
 \end{split}
\end{align}
Then, we can write
\begin{align}\label{eq: h_bar bound1}
\begin{split}
\| \overline{h}_t \| \leq& \|\Xi(t,T)\| \|\overline{h}_T\| \\
&+ \sum_{j=T-1}^{t} \|\Xi(t,j)\|  \|\left(I_d \otimes X_2(j)\right)\| \|\overline{m}_j\|.
 \end{split}
\end{align}
Now bound each term separately. For $t_2>t_1$, we have
\begin{align}\label{eq: transition matrix}
\begin{split}
 \Xi(t_1,t_2)&=\prod_{i=t_2-1}^{t_1} \Phi_i \otimes X_1(i)\\
 &= \left(\prod_{i=t_2-1}^{t_1} \Phi_i\right) \otimes \left(\prod_{i=t_2-1}^{t_1} X_1(i)\right).
\end{split}
\end{align}
Notice that, from \eqref{eq: h} and \eqref{eq: feedback gains}, we have $X_1(i)=A^{\T} + K_{i,x}^{\T} B^{\T}=A_c(i)^{\T}$. Then, we can write
\begin{align}\label{eq: transition matrix2}
\begin{split}
\prod_{i=t_2-1}^{t_1} X_1(i)= \prod_{i=t_2-1}^{t_1} {A_c(i)}^{\T} = \left(\prod_{i=t_1}^{t_2-1} A_c(i)\right)^{\T}.
\end{split}
\end{align}
Then, noting that $\|{\cdot}^{\T}\| = \|{\cdot}\|$, we can upper bound \eqref{eq: transition matrix}  as 
\begin{align}\label{eq: transition matrix bound}
\begin{split}
\left\| \Xi(t_1,t_2)\right\| &\leq \left(\prod_{i=t_2-1}^{t_1} \left\| \Phi_i\right\|\right)  \left(\left\|\prod_{i=t_1}^{t_2-1} A_c(i)\right\|\right) .
\end{split}
\end{align}
From \cite[Lemma B.1]{FC-GB-FP:22}, we have $\left\|\prod_{i=t_1}^{t_2-1} A_c(i)\right\|\leq \alpha\rho^{t_2-t_1}$, where $\alpha>0$ and $\rho\in (0,1)$ are constants. Further, using Assumption \ref{assump: linMDP}, we have for any $i\leq 0$
\begin{align}\label{eq: bound on phi_h}
\begin{split}
 \|\Phi_i\|_2&\leq {\|\Phi_i\|}_{\mathrm{F}} \\
 &=\sqrt{\trace{\left(\left(\Expectover{\mu_i}{\phi^{\T}(s_{i})}\right) \left(\Expectover{\mu_i}{\phi^{\T}(s_{i})}\right)^{\T}  \right)} }\\
 &=\sqrt{\sum_{j=1}^d \left(\Expectover{\mu_{j,i}}{\phi^{\T}(s_{i})}\right) \left(\Expectover{\mu_{j,i}}{\phi^{\T}(s_{i})}\right)^{\T}  } \\
  &=\sqrt{\sum_{j=1}^d {\|\Expectover{\mu_{j,i}}{\phi^{\T}(s_{i})}\|}_2^2  }\le	\sqrt{\sum_{j=1}^d \frac{1}{d}} =1 .
  \end{split}
\end{align}
Hence, we can re-write \eqref{eq: transition matrix bound} as
\begin{align}\label{eq: transition matrix bound2}
\left\| \Xi(t_1,t_2)\right\| &\leq \alpha \rho^{t_2-t_1}.
\end{align}
From Theorem \ref{thrm: Vfunction}, we have $\overline{h}_T=\Expectover{\mu_T}{h_T(s_{T})}=F \Expectover{\mu_T}{s_{T}}$. Then, we have
\begin{align}\label{eq: h_T bound}
\begin{split}
 \left\|\overline{h}_T\right\|\leq \left\|F\right\| \left\|\Expectover{\mu_T}{s_{T}}\right\| \overset{\text{(a)}}{\leq} \left\|F\right\|\Expectover{\mu_T}{\|s_{T}\|} \overset{\text{(b)}}{\leq} \left\|F\right\| \delta_s,
 \end{split}
\end{align}
where in step (a) we have used the Jensen's inequality, and in step (b) we have used Assumption \ref{assump: linMDP}. Next, we have
\begin{align}\label{eq: m_bar bound}
\begin{split}
\|\overline{m}_t\|=\|\Expectover{\mu_t}{s_t}\|&=\sqrt{\sum_{i=1}^d \left(\Expectover{\mu_{i,t}}{s_t}\right)^{\T} \left(\Expectover{\mu_{i,t}}{s_t}\right) }\\
&=\sqrt{\sum_{i=1}^d {\|\Expectover{\mu_{i,t}}{s_t}\|}^2 }\!\leq\! \sqrt{d \delta_s^2}\! =\! \delta_s \sqrt{d}.
\end{split}
 \end{align} 
 Let $\|X_2(t)\|\leq \overline{X}_2$ for $t\in \{0,\cdots, T-1\}$. Then, using \eqref{eq: transition matrix bound2}, \eqref{eq: h_T bound}, and \eqref{eq: m_bar bound} we can write \eqref{eq: h_bar bound1} as
\begin{align}\label{eq: h_bar bound2}
\begin{split}
 \|\overline{h}_t \| &\leq \left\|F\right\| \delta_s \alpha \rho^{T-t} + \overline{X}_2 \delta_s \alpha \sqrt{d}   \sum_{j=T-1}^t \rho^{j-t}\\
 &\overset{\text(c)}{=}\left\|F\right\| \delta_s \alpha \rho^{T-t} + \overline{X}_2 \delta_s \alpha \sqrt{d}   \sum_{k=0}^{T-t-1} \rho^{k}\\
 &=\left\|F\right\| \delta_s \alpha \rho^{T-t} + \overline{X}_2 \delta_s \alpha \sqrt{d}  \left(\frac{1-\rho^{T-t}}{1-\rho}\right)\\
 &\leq\left\|F\right\| \delta_s \alpha \rho^{T-t} +   \frac{\overline{X}_2 \delta_s \alpha \sqrt{d}}{1-\rho}.
 \end{split}
\end{align}
Now we bound $\overline{q}_t$ for $t\in \{0,\cdots ,T\}$. From \eqref{eq: q_bar dynamics}, we have
\begin{align}
\begin{split}
\overline{q}_t
&=\Phi_t \overline{q}_{t+1}+ v_t,
\end{split}
\end{align}
where,
\begin{align}\label{eq: vt}
\begin{split}
v_t&=
\Expectover{\mu_t}{s_t^{\T}Y_1(t) s_t}\\
&+ \begin{bmatrix}
 \overbar{h}_{t+1}^{\T} \left(\Expectover{\mu_{1,t}}{\phi(s_t)\phi(s_t)^{\T}} \otimes Y_2(t)\right)\overbar{h}_{t+1}\\
 \vdots \\
\overbar{h}_{t+1}^{\T} \left(\Expectover{\mu_{d,t}}{\phi(s_t)\phi(s_t)^{\T}} \otimes Y_2(t)\right)\overbar{h}_{t+1}
\end{bmatrix}\\
&+2 \Expectover{\mu_t}{\phi(s_t)^{\T} \otimes s_t^{\T}Y_3(t)} \overline{h}_{t+1}.
\end{split}
\end{align}
Given $\overline{q}_T$, we can write for $t\in \{0,\cdots,T-1\}$
\begin{align}\label{eq: q_bar at t}
 \overline{q}_t=\prod_{i=T-1}^{t} \Phi_i \overline{q}_T + \sum_{i=T-1}^t \prod_{j=i-1}^{t} \Phi_j v_i.
\end{align}
Then, we can upper bound $\|\overline{q}_t\|$ as
\begin{align}\label{eq: q_bar bound1}
\begin{split}
 \|\overline{q}_t\| &\leq \prod_{i=T-1}^{t} \|\Phi_i \|  \|\overline{q}_T\| + \sum_{i=T-1}^t \prod_{j=i-1}^{t} \|\Phi_j\| \|v_i\| \\
 & \overset{\text{(d)}}{\leq} \|\overline{q}_T\| + \sum_{i=T-1}^t  \|v_i\|,
\end{split}
\end{align}
where in step (d) we have used \eqref{eq: bound on phi_h}. Now we bound each term of $v_t$ for $t\!\in\! \{0,\cdots\! ,T\!-\!1\}$. We~start~with
\begin{align} 
\begin{split}
 \|\Expectover{\mu_{i,t}}{s^{\T}(t)Y_1(t) s_t}\| &\leq \Expectover{\mu_{i,t}}{\|s^{\T}(t)Y_1(t) s_t\|}\\
 &\leq  \Expectover{\mu_{i,t}}{{\|s_t\|}^2 \|Y_1(t)\| s_t\|} \\
 &\leq \delta_s^2 \|Y_1(t)\|,
\end{split}
\end{align}
for $i \in \{1, \cdots , d \}$. Then, we have 
\begin{align}\label{eq: v bound term 1} 
\begin{split}
 \|\Expectover{\mu_t}{s^{\T}(t)Y_1(t) s_t}\| &= \sqrt{\sum_{i=1}^d {\left\|\Expectover{\mu_{i,t}}{s^{\T}(t)Y_1(t) s_t}\right\|}^2}\\
 &\leq  \delta_s^2 \|Y_1(t)\| \sqrt{d}.
\end{split}
\end{align}
Next, for $i\in \{1, \cdots , d\}$, we have
\begin{align} 
\begin{split}
& \|\overline{h}_{t+1}^{\T} \left(\Expectover{\mu_{i,t}}{\phi(s_t) \phi^{\T}(s_t)} \otimes Y_2(t)\right) \overline{h}_{t+1} \|\\
 &\qquad\qquad\qquad\qquad\leq {\| \overline{h}_{t+1}\|}^2 {\|\phi(s_t)\|}^2 \|Y_2(t)\|\\
& \qquad\qquad\qquad\qquad\leq \frac{{\|\overline{h}_{t+1}\|}^2 \|Y_2(t)\|}{d}.
\end{split}
\end{align}
Then,
\begin{align}\label{eq: v bound term 2} 
\begin{split}
 &\left\| \begin{bmatrix}
 \overbar{h}_{t+1}^{\T} \left(\Expectover{\mu_{1,t}}{\phi(s_t)\phi(s_t)^{\T}} \otimes Y_2(t)\right)\overbar{h}_{t+1}\\
 \vdots \\
\overbar{h}_{t+1}^{\T} \left(\Expectover{\mu_{d,t}}{\phi(s_t)\phi(s_t)^{\T}} \otimes Y_2(t)\right)\overbar{h}_{t+1}
\end{bmatrix}\right\| \\
&\qquad=\sqrt{\sum_{i=1}^d {\left\|\overline{h}_{t+1}^{\T} \left(\Expectover{\mu_{i,t}}{\phi(s_t) \phi^{\T}(s_t)} \otimes Y_2(t)\right) \overline{h}_{t+1}\right\|}^2}\\
&\qquad \leq \frac{{\|\overline{h}_{t+1}\|}^2 \|Y_2(t)\|}{\sqrt{d}}.
\end{split}
\end{align}
Next, for $i\in \{1, \cdots , d\}$, we have
\begin{align}
\begin{split}
\left\|\Expectover{\mu_{i,t}\!\!}{\phi(s_t)^{\T} \otimes s_t^{\T}Y_3(t)} \right\|&\overset{\text{(e)}}{\leq} \Expectover{\mu_{i,t}\!\!}{\left\|\phi(s_t)^{\T} \otimes s_t^{\T}Y_3(t)\right\|}\\
&\leq \Expectover{\mu_{i,t}\!\!}{\left\|\phi(s_t)\right\| \|s_t\| \|Y_3(t)\| }\\
&\leq\frac{\delta_s \|Y_3(t)\|}{\sqrt{d}},
\end{split}
\end{align}
where in step (e) we have used Jensen's inequality. Then,
\begin{align}\label{eq: v bound term 3} 
\begin{split}
&\left\|\Expectover{\mu_t}{\phi(s_t)^{\T} \otimes s_t^{\T}Y_3(t)} \right\| \\
&\quad=\sqrt{\sum_{i=1}^d {\left\|\Expectover{\mu_{i,t}}{\phi(s_t)^{\T} \otimes s_t^{\T}Y_3(t)} \right\|}^2} \leq \delta_s \|Y_3(t)\|.
\end{split}
\end{align}
Let $\|Y_1(t)\|\leq \overline{Y}_1$, $\|Y_2(t)\|\leq \overline{Y}_2$, and $\|Y_3(t)\|\leq \overline{Y}_3$ for $t\in \{0,\cdots, T\}$. Using \eqref{eq: v bound term 1}, \eqref{eq: v bound term 2}, and \eqref{eq: v bound term 3}
\begin{align}\label{eq: v bound}
\begin{split}
 \|v_t\|\leq \delta_s^2 \overline{Y}_1 \sqrt{d} + \frac{\overline{Y}_2 {\|\overline{h}_{t+1}\|}^2}{\sqrt{d}} + \delta_s \overline{Y}_3 \|\overline{h}_{t+1}\|,
\end{split}
\end{align}
for $t\in \{0,\cdots, T\}$. From Theorem \ref{thrm: Vfunction}, we have $\overline{q}_T=\Expectover{\mu}{q_T(s_{T})}= \Expectover{\mu}{{s_{T}}^{\T} M s_{T}}$. Then, $\|\overline{q}_T\|\leq \delta_s^2 \|M\|$. Then, we can re-write \eqref{eq: q_bar bound1} as
\begin{align}\label{eq: q_bar bound2}
 \begin{split}
\|\overline{q}_t\|\leq \delta_s^2 \|M\| + \delta_s^2 \overline{Y}_1 \sqrt{d} + \frac{\overline{Y}_2 {\|\overline{h}_{t+1}\|}^2}{\sqrt{d}} + \delta_s \overline{Y}_3 \|\overline{h}_{t+1}\|.
\end{split}
\end{align}
\end{proof}
\end{theorem}
%
%
%
Following the same notation as \eqref{eq: theta}, let the true parameter be denoted by $\theta^*$, which is written as
\begin{align}\label{eq: theta_star}
\theta^*_t= \begin{bmatrix}
{ \theta^*_{1,t}}^{\T}&
 \cdots&
{\theta^*_{d,t}}^{\T}
\end{bmatrix}^{\T},
\quad
\text{where}
\quad
\theta^*_{i,t} =\begin{bmatrix}
 \overline{h}_{i,t}\\
 \overline{q}_{i,t}
\end{bmatrix}.
\end{align}
where $\overline{h}_{i,t}\in \mathbb{R}^n$ and $\overline{q}_{i,t}\in \mathbb{R}$ are the components of $\overline{h}_{t}$ and $\overline{q}_{t}$ in Theorem \ref{thrm: true parameters} for $i\in\{1,\cdots,d\}$ and $t\in\{0,\cdots,T\}$.
\begin{corollary}{\bf \emph{(bound on $\theta^*_t$)}}\label{cor: theta_star bound}
Let $\theta^*_t$ be as in \eqref{eq: theta_star}. Then, under the same assumptions of Theorem~\ref{thrm: true parameters} and Theorem~\ref{thrm: true parameters bounds}, we have $\|\theta_t^*\|\leq c_{\theta}\sqrt{d}$ for $t\in\{0,\cdots,T\}$, where $c_{\theta}>0$ is independent of $d$.
\begin{proof}
Theorem \ref{thrm: true parameters bounds} implies that for $t\in\{0,\cdots,T\}$,
\begin{align}\label{eq: true parameters bounds 2}
\|\overline{h}_t\| \leq a_h+ b_h\sqrt{d},\quad \|\|\overline{q}_t\| \leq a_q+ b_q\sqrt{d},
\end{align}
where $a_h>0$, $b_h>0$, $a_q>0$, and $b_q>0$ are independent of $d$. Then, we can bound $\theta^*_t$ in \eqref{eq: theta_star} as
\begin{align}
\begin{split}
\| \theta_t^*\|=&\sqrt{\sum_{i=1}^d ({\theta_{i,t}^*})^{\T}\theta_{i,t}^* }=\sqrt{\sum_{i=1}^d ({\overline{h}_{i,t}})^{\T}\overline{h}_{i,t}+({\overline{q}_{i,t}})^{\T}\overline{q}_{i,t} }\\
=&\sqrt{(\overline{h}_t)^{\T}\overline{h}_t + (\overline{q}_t)^{\T}\overline{q}_t}=\sqrt{\|\overline{h}_t\|^2+\|\overline{q}_t\|^2}\\
\leq&\|\overline{h}_t\|+\|\overline{q}_t\|\leq \underbrace{(a_h+b_h+a_q+b_q)}_{c_{\theta}}\sqrt{d}.
\end{split}
\end{align}
\end{proof}
\end{corollary}
Corollary \ref{cor: theta_star bound} implies that choosing the projection radius in Algorithm \ref{alg: alg1} as $R_{\theta}\geq c_{\theta}\sqrt{d}$ guarantees that $\theta^*_t$ belongs to the projection ball for all $t$.

\section{Proof of Theorem \ref{thrm: ISS}} \label{app: ISS}
Let $A_c(t)=A+B K_{t,x}$ and let $\varphi(t_2,t_1)=\prod_{i=t_1}^{t_2-1}A_c(i)$ denote the state transition matrix from $t_1$ to $t_2$.\footnote{The matrix multiplication is performed from the left, i.e., $A(t_1)$ appears as the rightmost matrix in the product.} Let $\pi_t^{\ell}(x_t,s_t)=K_{t,x} x^{\ell}_t + K_{t,s} s^{\ell}_t + K_{t,h} \left({\phi\left(s^{\ell}_t\right)}^{\T} \otimes Z \right)\theta_{t+1}^{\ell}$ denote the policy learned from Algorithm \ref{alg: alg1} at episode $\ell$ and time $t$, where $Z=[I_n, 0_{n\times1}]$. Then, the evolution of $x_t$ in system \eqref{eq: system} under the policy $\{\pi_1^{\ell},\cdots\! , \pi_{t-1}^{\ell}\}$ for $t\in\{0,\cdots\!,T\}$ is written as 
\begin{align}\label{eq: evolution of x}
 x^{\ell}_t=\varphi(t,0)x^{\ell}_0 + \sum_{i=0}^{t-1} \varphi(t,i+1) B \overline{u}^{\ell}_i, 
\end{align}
where $x^{\ell}_0$ is the initial state at episode $\ell$ and $\overline{u}^{\ell}_i=K_{i,s}s^{\ell}_i + K_{i,h}\left(\phi(s^{\ell}_i)\otimes Z\right)\theta_{i+1}^{\ell}$. Then,
\begin{align}\label{eq: x bound}
\begin{split}
 \|x^{\ell}_t\|\leq &\|\varphi(t,0)\| \|x^{\ell}(0)\|\\
 & + \|B\| \sum_{i=0}^{t-1} \|\varphi(t,i+1)\| \left(\sup_{0\leq j\leq t-1} \|\overline{u}^{\ell}_j\|\right).
 \end{split}
\end{align}
From \cite[Lemma B.1]{FC-GB-FP:22}, we have $\|\varphi(t_2,t_1)\|\leq \alpha\rho^{t_2-t_1}$ where $\alpha>0$ and $\rho\in (0,1)$ are constants. Then, we can write \eqref{eq: x bound} as
\begin{align}\label{eq: x bound 2}
\begin{split}
  \|x^{\ell}_t&\|\leq \alpha \rho^t \|x^{\ell}_0\| + \alpha \|B\| \sum_{i=0}^{t-1} \rho^{t-i-1} \underbrace{\left(\sup_{0\leq j\leq t-1} \|\overline{u}^{\ell}_j\|\right)}_{u_{\infty}^{\ell}}\\
  &\overset{\text{(a)}}{=} \alpha \rho^t \|x^{\ell}_0\| + \alpha \|B\| \sum_{k=0}^{t-1} \rho^{k} u_{\infty}^{\ell}\\
  &= \alpha \rho^t \|x^{\ell}_0\| + \alpha \|B\| \left(\frac{1-\rho^{t-1}}{1-\rho}\right) u_{\infty}^{\ell}\\
  &\leq \alpha \rho^t \|x^{\ell}_0\| + \alpha \|B\| \left(\frac{1}{1-\rho}\right) u_{\infty}^{\ell},
  \end{split}
\end{align}
where in step (a), we have changed the index in the sum to $k=t-i-1$. Next, we bound on $u_{\infty}^{\ell}$. Let $\|\theta_t^{\ell}\|\leq R_{\theta}$, $\|K_{t,s}\|\leq \overline{K}_s$, and $\|K_{t,h}\|\leq \overline{K}_h$ for $t\in\{0, \cdots , T - 1\}$ and episode $\ell$. Then we have,
\begin{align}
\begin{split}
\|\overline{u}^{\ell}_t\| & \leq \|K_{t,s}\| \|s^{\ell}_t\| + \|K_{t,h}\|\left(\phi(s^{\ell}_t)\otimes Z\right)\| \|\theta_{t+1}^{\ell}\|\\
&\leq \overline{K}_s \delta_s + \overline{K}_h \|\phi(s^{\ell}_t)\| \|Z\| R_{\theta}\\
&\leq \overline{K}_s \delta_s + \frac{\overline{K}_h  R_{\theta}}{\sqrt{d}}.
\end{split}
\end{align}
Since the above bound is uniform for $t\in\{0,\cdots , T-1\}$, we have $u_{\infty}^{\ell} \leq \overline{K}_s \delta_s + \frac{\overline{K}_h  R_{\theta}}{\sqrt{d}}$. The proof follows by substituting the bound of $u_{\infty}^{\ell}$ in \eqref{eq: x bound 2}.

\section{Proof of Theorem \ref{thrm: regret bound}} \label{app: regret bound}
We begin by presenting the following technical Lemmas.
\begin{lemma}{\bf \emph{(Chernoff bound)}}\label{lemma: norm_Lambda_inv}
Let $X_t\!=\!\sum_{i=1}^{t-1} {z_i}{z_i}^{\T} \!+\!\gamma I_{p}$, where $z_i\! \in\! \mathbb{R}^{p}$~and $\gamma\!>\!0$. Let $\expect{zz^{\T}}\!\succeq\! \alpha I_p$ with $\alpha>0$, and $\|z\|\leq \zeta$. Let $\delta \in [0,1]$ and assume $t\geq (8\zeta^2 \log(p/\delta))/\alpha$. Then, with probability at least $1-\delta$, the minimum eigenvalue of $X_t$~satisfies
\begin{align*}
\lambda_{\text{min}}\left(X_t\right)\geq \gamma + \frac{(t-1)\alpha}{2}.
\end{align*}
\begin{proof}
Let $\alpha = \lambda_{\text{min}}(\expect{zz^{\T}})$. Define 
\begin{align*}
 \mu_{\text{min}}\triangleq &\lambda_{\text{min}}\left(\sum_{i=1}^{t-1}\expect{zz^{\T}}\right)=\lambda_{\text{min}}\left((t-1)\expect{zz^{\T}}\right)\\
 =&(t-1)\lambda_{\text{min}}\left(\expect{zz^{\T}}\right)=(t-1)\alpha.
\end{align*}
Further, we have $z_i z_i^{\T}\succeq 0$ and $\lambda_{\text{max}}\left(z_i z_i^{\T}\right)\!=\!\|z_i\|^2\!\leq \!\zeta^2$. Then, using \cite[Theorem 1.1]{JAT:12}, we have
\begin{align}\label{eq: chernoff bound}
\begin{split}
 &\mathbb{P} \left(\lambda_{\text{min}}\left(\sum_{i=1}^{t-1}z_iz_i^{\T}\right)\leq (1-\varepsilon) (t-1)\alpha \right)\\
 &\qquad \qquad \qquad \qquad\qquad\leq p\left(\frac{\exp(-\varepsilon)}{(1-\varepsilon)^{1-\varepsilon}}\right)^{\frac{(t-1)\alpha}{\zeta^2}},\\
 &\qquad \qquad \qquad \qquad\qquad \overset{\text{(a)}}{\leq} p\exp\left(\frac{-\varepsilon^2(t-1)\alpha}{2\zeta^2}\right),
 \end{split}
\end{align}
for $\varepsilon \in [0,1]$, where in step (a) we have used $\frac{\exp(-\varepsilon)}{(1-\varepsilon)^{1-\varepsilon}}\leq \exp\left(-\varepsilon^2/2\right)$ for $\varepsilon \in (0,1)$. Choose $\varepsilon=0.5$, then we write \eqref{eq: chernoff bound} as
\begin{align}\label{eq: chernoff bound 2}
\begin{split}
 &\mathbb{P} \left(\lambda_{\text{min}}\left(\sum_{i=1}^{t-1}z_iz_i^{\T}\right)\leq \frac{(t-1)\alpha}{2} \right)\\
 &\qquad \qquad \qquad \qquad\qquad\leq p\exp\left(\frac{-(t-1)\alpha}{8\zeta^2}\right).
 \end{split}
\end{align}
Let $p\exp\left(\frac{-(t-1)\alpha}{8\zeta^2}\right) \leq \delta$, then we have
\begin{align*}
 t\geq \frac{8 \zeta^2 \log{(p/\delta)}}{\alpha}.
\end{align*}
Then, with probability at least $1-\delta$ we have
\begin{align*}
 \lambda_{\text{min}}\left(\sum_{i=1}^{t-1}z_iz_i^{\T}\right)\geq \frac{(t-1)\alpha}{2}.
\end{align*}
Finally, we have
\begin{align*}
 \lambda_{\text{min}}\left(X_t\right)\geq &\lambda_{\text{min}}\left(\sum_{i=1}^{t-1}z_iz_i^{\T}\right) +\gamma \geq \frac{(t-1)\alpha}{2} +\gamma.
\end{align*}
\end{proof}
\end{lemma}
\begin{lemma}\label{lemma: expected_psi_psiT}
Consider the system \eqref{eq: system} and the Markov process \eqref{eq: environment}. Let Assumption \ref{assump: linMDP} be satisfied, and let 
\begin{align*}
 x_{t+1}=\varphi(t+1,0)x_0 + \sum_{i=0}^{t} \varphi(t+1,i+1) B \overline{u}_i(s_i), 
\end{align*}
with $x_0\sim \mathcal{N}(0,\Sigma_0)$, and $\overline{u}_i(s_i)$ is an arbitrary input that depends on $s_i$ and is independent of $x_0$. Let 
\begin{align*}
 \psi_t=\phi(s_t)\otimes 
\begin{bmatrix}
 2x_{t+1}\\1
\end{bmatrix}.
\end{align*}
Assume $\Sigma_0\succ 0$ and $\varphi(t+1,0)$ is nonsingular for $t\in\{0,\cdots, T-1\}$. Then, $\expect{\psi_t \psi_t^{\T}} \succ 0$ for $t\in\{0,\cdots, T-1\}$.
\begin{proof}
We begin by writing 
\begin{align*}
 \begin{split}
\psi_t= \phi(s_t)\otimes
\begin{bmatrix}
 2 x_{t+1}\\ 1
\end{bmatrix} 
=
\begin{bmatrix}
2 \phi(s_t)\otimes  x_{t+1}\\ 
\phi(s_t)
\end{bmatrix}.
\end{split}
\end{align*}
Then,
\begin{align*}
 \begin{split}
&\psi_t {\psi_t}^{\T}\\
&=
\begin{bmatrix}
4 \phi(s_t)\phi(s_t)^{\T} \!\otimes\!  x_{t+1} {x_{t+1}}^{\T} & 2\phi(s_t)\phi(s_t)^{\T} \!\otimes\!  x_{t+1}\\
2\phi(s_t)\phi(s_t)^{\T} \!\otimes\!  {x_{t+1}}^{\T} & \phi(s_t){\phi(s_t)}^{\T}
\end{bmatrix}.
\end{split}
\end{align*}
Taking the expectation, we get
\begin{align*}
 \begin{split}
&\expect{\psi_t {\psi_t}^{\T}}\\
&\!\!=\!\!
\left[\begin{smallarray}{cc}
4 \expect{\phi(s_t)\phi(s_t)^{\T} \!\otimes \!x_{t+1} {x_{t+1}}^{\T}} & ~~2\expect{\phi(s_t)\phi(s_t)^{\T} \!\otimes \!x_{t+1}}\\
2\expect{\phi(s_t)\phi(s_t)^{\T} \otimes  {x_{t+1}}^{\T}} & \expect{\phi(s_t){\phi(s_t)}^{\T}}
\end{smallarray}\right]. 
\end{split}
\end{align*}
From Assumption \ref{assump: linMDP}, we have $\expect{\phi(s_t){\phi(s_t)}^{\T}}\succ 0$ for all $t$. For notational convenience we denote $\Sigma_{\phi}=\expect{\phi(s_t){\phi(s_t)}^{\T}}$. We apply the Schur complement
\begin{align}\label{eq: schur_comp}
\begin{split}
 S&= 4\expect{\phi(s_t)\phi(s_t)^{\T} \otimes  x_{t+1} {x_{t+1}}^{\T}}\\
  &-\! 4\expect{\phi(s_t)\phi(s_t)^{\T} \!\otimes \! x_{t+1}} \Sigma_{\phi}^{-1} \expect{\phi(s_t)\phi(s_t)^{\T}\! \otimes\!  {x_{t+1}}^{\T}}\!.
\end{split}
\end{align}
Showing $\expect{\psi_t {\psi_t}^{\T}}\succ 0$ boils down to showing that $S\succ 0$. From \eqref{eq: schur_comp}, we have
\begin{align}\label{eq: schur_comp_1st}
\begin{split}
 &\expect{\phi(s_t)\phi(s_t)^{\T} \otimes  x_{t+1} {x_{t+1}}^{\T}}\\
&\qquad\qquad \overset{\text{(a)}}{=}\expect{\expect{\phi(s_t)\phi(s_t)^{\T} \otimes  x_{t+1} {x_{t+1}}^{\T}\big| s_t}}\\
&\qquad\qquad=\expect{\phi(s_t)\phi(s_t)^{\T} \otimes  \expect{x_{t+1} {x_{t+1}}^{\T} \big| s_t}}\\
&\qquad\qquad=\expect{\phi(s_t){\phi(s_t)}^{\T} \otimes  \Sigma_{x|s}} \\
&\qquad\qquad\quad+ \expect{\phi(s_t){\phi(s_t)}^{\T} \otimes  \mu_{x}(s_t) {\mu_{x}(s_t)}^{\T}},
\end{split}
\end{align}
where in step (a) we used the law of total expectation,~and 
\begin{align*}
 \begin{split}
 \Sigma_{x|s}\!=&\expect{\left(x_{t+1}\!-\!\expect{x_{t+1}|s_t}\right) \left(x_{t+1}\!-\!\expect{x_{t+1}|s_t}\right)^{\T}\big|s_t}\!,\\
 \mu_{x}(s_t)=&\expect{x_{t+1}|s_t}.
\end{split}
\end{align*}
For notational convenience, let $z=\phi(s_t)\otimes \mu_x(s_t)$. Substituting \eqref{eq: schur_comp_1st} in \eqref{eq: schur_comp}, we get
\begin{align}\label{eq: schur_comp_cont}
\begin{split}
 S=& \underbrace{4\expect{\phi(s_t){\phi(s_t)}^{\T} \otimes  \Sigma_{x|s} }}_{S_1} \\
 &+ \underbrace{4\expect{z z^{\T}} - 4 \expect{z {\phi(s_t)}^{\T} }{\Sigma_{\phi}}^{-1} \expect{\phi(s_t) z^{\T}}}_{S_2}.
\end{split}
\end{align}
We have $S_1\succeq 0$ since $\phi(s_t)\phi(s_t)^{\T} \succeq 0$ and $\Sigma_{x|s}\succeq 0$. From \eqref{eq: schur_comp_cont}, we have
\begin{align*}
\begin{split}
S_2=&4\expect{z z^{\T}} - 4 \expect{z {\phi(s_t)}^{\T} }{\Sigma_{\phi}}^{-1} \expect{\phi(s_t) z^{\T}}\\
\overset{\text{(b)}}{=}&4\expect{z z^{\T}} - 4 \expect{z {\phi(s_t)}^{\T}{\Sigma_{\phi}}^{-1}\expect{\phi(s_t) z^{\T}}}\\
&+4\expect{\expect{z {\phi(s_t)}^{\T}}{\Sigma_{\phi}}^{-1}\phi(s_t) z^{\T}} \\
&-4\expect{\expect{z {\phi(s_t)}^{\T}}{\Sigma_{\phi}}^{-1}\phi(s_t) z^{\T}}\\
\overset{\text{(c)}}{=}&4\expect{z z^{\T}} - 4 \expect{z {\phi(s_t)}^{\T}{\Sigma_{\phi}}^{-1}\expect{\phi(s_t) z^{\T}}}\\
&+4\expect{z {\phi(s_t)}^{\T}}{\Sigma_{\phi}}^{-1}\Sigma_{\phi}{\Sigma_{\phi}}^{-1}\expect{\phi(s_t) z^{\T}} \\
&-4\expect{\expect{z {\phi(s_t)}^{\T}}{\Sigma_{\phi}}^{-1}\phi(s_t) z^{\T}}\\
=&4 \mathbb{E}\Big[z z^{\T} - z {\phi(s_t)}^{\T}{\Sigma_{\phi}}^{-1}\expect{\phi(s_t) z^{\T}} \\
&- \expect{z {\phi(s_t)}^{\T}}{\Sigma_{\phi}}^{-1}\phi(s_t) z^{\T}\\
& +\expect{z {\phi(s_t)}^{\T}}{\Sigma_{\phi}}^{-1}\phi(s_t){\phi(s_t)}^{\T} {\Sigma_{\phi}}^{-1}\expect{\phi(s_t) z^{\T}} \Big]\\
=&4 \mathbb{E}\Big[\left(z- \expect{z {\phi(s_t)}^{\T}}{\Sigma_{\phi}}^{-1}\phi(s_t)  \right)\\
&\qquad \qquad \qquad . \left(z- \expect{z {\phi(s_t)}^{\T}}{\Sigma_{\phi}}^{-1}\phi(s_t)\right)^{\T} \Big],
\end{split}
\end{align*}
where in step (b) we have added and subtracted the term $4\expect{\expect{z {\phi(s_t)}^{\T}}{\Sigma_{\phi}}^{-1}\phi(s_t) z^{\T}}$, and in step (c) we have used $I= \Sigma_{\phi} {\Sigma_{\phi}}^{-1}$. Then, we have $S=S_1+S_2\succeq 0$. From \eqref{eq: evolution of x} we have $x_{t+1}=\varphi(t+1,0) x_0 +\sum_{i=0}^{t} \varphi(t+1,i+1) B \overline{u}_i$, hence, $\expect{x_{t+1}|s_t} =  \expect{ \sum_{i=0}^{t} \varphi(t+1,i+1) B \overline{u}_i | s_t}$. For notational convenience, let $\widetilde{u}(t)=\sum_{i=0}^{t} \varphi(t+1,i+1) B \overline{u}_i$. Then we get
\begin{align*}
 \Sigma_{x|s}&=\varphi(t+1,0) \Sigma_0 {\varphi(t+1,0)}^{\T}\\
 &+ \expect{\left(\widetilde{u}(t) - \expect{\widetilde{u}(t) | s_t} \right) \left(\widetilde{u}(t) - \expect{\widetilde{u}(t) | s_t} \right)^{\T} \big| s_t}\\
&\succeq \varphi(t+1,0) \Sigma_0 {\varphi(t+1,0)}^{\T}. 
\end{align*}
Since $\Sigma_0 \succ 0$ and $\varphi(t+1,0)$ is nonsingular for all $t$, then, $\varphi(t+1,0) \Sigma_0 {\varphi(t+1,0)}^{\T}\succ 0$. Then, we have $S_1=4\expect{\phi(s_t){\phi(s_t)}^{\T} \otimes  \Sigma_{x|s} }\succeq 4\expect{\phi(s_t){\phi(s_t)}^{\T} } \otimes \left(\varphi(t+1,0) \Sigma_0 {\varphi(t+1,0)}^{\T}\right) \succ 0$ since $\expect{\phi(s_t){\phi(s_t)}^{\T} }~\succ 0$, and $\varphi(t+1,0) \Sigma_0 {\varphi(t+1,0)}^{\T}$ is independent of $s$. Therefore, $S\succ 0$, which implies $\expect{\psi_t {\psi_t}^{\T}}\succ~0$ for $t \in \{0,\cdots,T-1\}$.
\end{proof}
\end{lemma}
%
%
%
%
%
%
%
%
Now we present the proof of Theorem \ref{thrm: regret bound}. For notational convenience, we denote $\phi\left(s^i_t\right)$ and $Y\left(x^{i}_t,u^{i}_t\right)$ by $\phi_t^i$ and $Y_t^i$, respectively. We have from the expression of $\theta_{t+1}^{\ell}$ in \eqref{eq: LSVI sol}
\begin{align*}
 \epsilon_{t+1}^{\ell}(x_{t+1}^i,s_{t+1}^i)&=2 x_{t+1}^i h_{t+1}^{\ell}(s_{t+1}^i) +q_{t+1}^{\ell}(s_{t+1}^i)\\
 &= 
\underbrace{\begin{bmatrix}
 2x_{t+1}^i & 1
\end{bmatrix}}_{(y_{t}^i)^{\T}}
\underbrace{\begin{bmatrix}
 h_{t+1}^{\ell}(s_{t+1}^i)\\
 q_{t+1}^{\ell}(s_{t+1}^i)
\end{bmatrix}}_{v_{t+1}^{\ell}(s_{t+1}^i)}.
\end{align*}
We can derive an upper bound on $|\epsilon_{t+1}^{\ell}(x_{t+1}^i,s_{t+1}^i)|$ as
\begin{align}\label{eq: epsilon bound}
\begin{split}
 |\epsilon_{t+1}^{\ell}(x_{t+1}^i,s_{t+1}^i)| \leq &2\|x_{t+1}^i\| \|h_{t+1}^{\ell}(s_{t+1}^i)\|\\
 &+ \|q_{t+1}^{\ell}(s_{t+1}^i)\|.
\end{split} 
\end{align}
Using \eqref{eq: h in terms of theta} we can write
\begin{align}\label{eq: h bound}
\begin{split}
 \|h^{\ell}_{t+1}(s_{t+1}^i)\|\leq& \|X_1(t+1)\| \|\phi_{t+1}^i\| \|\theta_{t+2}^{\ell}\|\\
 & + \|X_2(t+1)\| \|s_{t+1}^i\| \\
 &\leq  \frac{\overline{X}_1  R_{\theta}}{\sqrt{d}} + \overline{X}_2 \delta_s,
 \end{split}
 \end{align}
where $\|X_1(t)\|\leq \overline{X}_1$ and $\|X_2(t)\|\leq \overline{X}_2$ for all $t\in \{0,\cdots,T\}$. From \eqref{eq: q in terms of theta} we can write
 \begin{align}\label{eq: q bound}
\begin{split}
 q^{\ell}_{t+1}(s_{t+1}^i) \leq& \|\phi^i_{t+1}\| \|\theta_{t+2}^{\ell}\| + \| s_{t+1}^i\|^2 \|Y_1(t+1)\|\\
 &+ \|\theta^{\ell}_{t+2}\|^2 \|\phi^i_{t+1}\|^2 \|Y_2(t+1)\| \\
 &+2 \|s^i_{t+1}\| \|Y_3(t+1)\| \|\phi^i_{t+1}\| \|\theta^{\ell}_{t+2}\| \\
 \leq& \frac{R_{\theta}}{\sqrt{d}} + \delta_s^2 \overline{Y}_1+ \frac{R_{\theta}^2 \overline{Y}_2}{d} +2 \frac{\delta_s \overline{Y}_3 R_{\theta}}{\sqrt{d}},
 \end{split}
\end{align}
where $\|Y_1(t)\|\leq \overline{Y}_1$, $\|Y_2(t)\|\leq \overline{Y}_2$, and $\|Y_3(t)\|\leq \overline{Y}_3$ for all $t\in \{0,\cdots,T\}$. From Theorem \ref{thrm: ISS}, we have
\begin{align*}
  ||x^{\ell}_t||\leq  \alpha\rho^t \|x^{\ell}(0)\| + \frac{\alpha\|B\|}{1-\rho}\left(\overline{K}_s\delta_s + \frac{\overline{K}_hR_{\theta}}{\sqrt{d}}\right),
\end{align*}
 for $t\in \{0,\cdots, T\}$ and $\ell \in \{1,\cdots, L\}$, with $\alpha>0$ and $ 0<\rho<1$. Define 
\begin{align}\label{eq: x bound}
 \overline{x} \!=\! \!\sup_{\substack{t\in \{0,\cdots, T\}\\ \ell \in \{1,\cdots, L\}}} \left\{ \alpha\rho^t \|x^{\ell}(0)\| + \frac{\alpha\|B\|}{1-\rho}\left(\overline{K}_s\delta_s + \frac{\overline{K}_hR_{\theta}}{\sqrt{d}}\right) \right\}\!.
\end{align}
Substituting \eqref{eq: h bound}, \eqref{eq: q bound}, and \eqref{eq: x bound} in \eqref{eq: epsilon bound}, we get
\begin{align}\label{eq: epsilon bound 2}
\begin{split}
|\epsilon_{t+1}^{\ell}(x_{t+1}^i,s_{t+1}^i)| \leq & \left(\frac{2\overline{x} \overline{X}_1 +1 + 2\delta_s \overline{Y}_3  } {\sqrt{d}} \right) R_{\theta}\\
& + \frac{\overline{Y}_2}{d}R_{\theta}^2 + 2 \overline{X}_2 \overline{x} \delta_s + \overline{Y}_1\delta_s^2.
\end{split}
\end{align}
Further, we can bound 
\begin{align}\label{eq: bound on psi}
\begin{split}
 \| \psi_{t}^i\|=\|(Y_t^i)^{\T} \phi_t^i\|&\leq \left\|
\begin{bmatrix}
2x_{t+1}^i & 1 
\end{bmatrix}\right\| \|\phi_t^i\|\\
&\leq\sqrt{\frac{4 \overline{x}^2 +1} {d}} \triangleq \delta_{\psi}.
\end{split}
\end{align}
Next, using the expression of $\theta_{t+1}^{\ell}$ in \eqref{eq: LSVI sol}, we write
\begin{align}\label{eq: delta_theta1}
\theta^{\ell}_{t+1}\!-\!\theta^*_{t+1}\!=\! (\Lambda_t^{\ell})^{-1} \sum_{i=1}^{\ell-1}\left(Y_t^i\right)^{\T} \phi_t^i \epsilon_{t+1}^{\ell}(x_{t+1}^i,s_{t+1}^i)\! -\!\theta^*_{t+1},
\end{align}
We re-write \eqref{eq: delta_theta1} as 
\begin{align}\label{eq: delta_theta2}
\begin{split}
  \theta^{\ell}_{t+1}\!-\!\theta^*_{t+1}\!=& (\Lambda_t^{\ell})^{-1} \Big(\sum_{i=1}^{\ell-1}(Y_t^i)^{\T} (\phi_t^i) (y_{t}^i)^{\T}  v_{t+1}^{\ell}(s_{t+1}^i)\\
  &\qquad \qquad \qquad \qquad -\Lambda_{t}^{\ell}\theta_{t+1}\Big)\\
  =&(\Lambda_t^{\ell})^{-1} \Bigg(\sum_{i=1}^{\ell-1}\left(Y_t^i\right)^{\T} (\phi_t^i) \left(y_{t}^i\right)^{\T}   v_{t+1}^{\ell}(s_{t+1}^i) \\
  & \quad\qquad- \sum_{i=1}^{\ell-1}(Y_t^i)^{\T} (\phi_t^i)(\phi_t^i)^{\T} Y_t^i \theta_{t+1}\Bigg)\\
  & - \lambda  (\Lambda_t^{\ell})^{-1} \theta^*_{t+1}.
  \end{split}
\end{align}
From \eqref{eq: delta_theta2}, we expand the term as
\begin{align}\label{eq: delta_theta3}
\begin{split}
 \left(\phi_t^i\right)^{\T} Y_t^i \theta^*_{t+1}=&\sum_{j=1}^d \phi_{j,t}^i\left(y_t^i\right)^{\T} 
\begin{bmatrix}
 \overline{h}_{j,t+1}^*\\
 \overline{q}_{j,t+1}^*
\end{bmatrix}\\
=&\left(y^i_{t}\right)^{\T}\sum_{j=1}^d \phi_{j,t}^i
\Expectover{\mu_j}{\begin{bmatrix}
h_{t+1}^*\left(s_{t+1}\right)\\
q_{t+1}^*\left(s_{t+1}\right)
\end{bmatrix}}\\
=&\left(y^i_{t}\right)^{\T}\sum_{j=1}^d \phi_{j,t}^i
\int_{\mathcal{S}}
 \underbrace{\begin{bmatrix}
h^*_{t+1}\left(s_{t+1}\right)\\
q^*_{t+1}\left(s_{t+1}\right)
\end{bmatrix}}_{v^*_{t+1}(s_{t+1})}
\mu_j\left(d s_{t+1}\right)\\
=&\left(y^i_{t}\right)^{\T}\int_{\mathcal{S}} v^*_{t+1}(s_{t+1}) \sum_{j=1}^d \phi_{j,t}^i\mu_j\left(d s_{t+1}\right)\\
=&\left(y^i_{t}\right)^{\T}\Expect{}{v^*_{t+1}(s_{t+1})}{s^i_t}.
\end{split}
\end{align}
For notational convenience, we use $\psi_t^i=(Y_t^i)^{\T}\phi^i_t$. Substituting \eqref{eq: delta_theta3} in \eqref{eq: delta_theta2}, we get
\begin{align}\label{eq: delta_theta4}
\begin{split}
  &\theta^{\ell}_{t+1}-\theta^*_{t+1} \\
  &=\!(\Lambda_t^{\ell})^{-1} \Bigg(\sum_{i=1}^{\ell-1}\psi_t^i(y_{t}^i)^{\T} \Big( v_{t+1}^{\ell}(s_{t+1}^i)\\
  &\qquad \qquad \qquad  - \Expect{}{v^*_{t+1}(s_{t+1})}{s^i_t}\Big)\Bigg)- \lambda  (\Lambda_t^{\ell})^{-1} \theta^*_{t+1}\\
  &=\!\underbrace{(\Lambda_t^{\ell})^{-1}  \Big(\sum_{i=1}^{\ell-1}\psi_t^i (y_{t}^i)^{\T} \big( v_{t+1}^{\ell}(s_{t+1}^i) \!-\! \Expect{}{v_{t+1}^{\ell}(s_{t+1})}{s^i_t}\big)\Big)}_{r_1}\\
  &+\!\underbrace{(\Lambda_t^{\ell})^{-1}  \Big(\sum_{i=1}^{\ell-1}\psi_t^i  (y_{t}^i)^{\T} \big( \Expect{}{v_{t+1}^{\ell}(s_{t+1})\! -\!v_{t+1}^{*}(s_{t+1})}{s^i_t}\big)\Big)}_{r_2}\\
  &\qquad-\underbrace{ \lambda  (\Lambda_t^{\ell})^{-1} \theta^*_{t+1}}_{r_3} .
  \end{split}
\end{align}
Let $\xi_{i,t+1}^{\ell}=(y_t^i)^{\T}v_{t+1}^{\ell}(s_{t+1}^i) - \Expect{}{(y_t^i)^{\T}v_{t+1}^{\ell}(s_{t+1})}{s^i_t}$. We have $\expect{\xi_{i,t+1}^{\ell} | \mathcal{F}_t^{\ell-1}}=0$. Since $(y_t^i)^{\T}v_{t+1}^{\ell}(s_{t+1}^i)$ is bounded (see \eqref{eq: epsilon bound 2}), then $\xi$ is $\sigma$-subGaussian with 
\begin{align}\label{eq: sigma bound}
\begin{split}
 \sigma =&\left(\frac{2\overline{x} \overline{X}_1 +1 + 2\delta_s \overline{Y}_3  } {\sqrt{d}} \right) R_{\theta}\\
& + \frac{\overline{Y}_2}{d}R_{\theta}^2 + 2 \overline{X}_2 \overline{x} \delta_s + \overline{Y}_1\delta_s^2.
\end{split}
\end{align}
Then, using \cite[Theorem~1]{YAY-DP-CS:11}, we have with probability at least $1-\delta$, we have
\begin{align}\label{eq: bounding r1}
\begin{split}
 &{\left\|\sum_{i=1}^{\ell-1}(\psi_t^i (y_{t}^i)^{\T} \left( v_{t+1}^{\ell}(s_{t+1}^i) - \Expect{}{ v_{t+1}^{\ell}(s_{t+1})}{s^i_t}\right)\right\|}^2_{\left({\Lambda_t^{\ell}}\right)^{-1}}\\
 &\qquad \leq 2 \sigma^2 \left(\log{\left( \sqrt{\frac{\det{\left(\Lambda_t^{\ell} \right)}}{\det{\left(\Lambda_t^{1} \right)}}}\right)} +\log{\left(\frac{1}{\delta} \right)}\right).
\end{split}
\end{align}
Recall from Alg. \ref{alg: alg1}, we have 
\begin{align*}
\begin{split}
 \Lambda_t^1&= \lambda I_{d(n+1)},\\
 \Lambda_t^{\ell} &= \sum_{i=1}^{\ell-1}{\left(Y^i_t\right)}^{\T} \phi^i_t {\left(\phi_t^i\right)}^{\T} Y_t^i + \lambda I_{d(n+1)}\\
 &=\lambda\underbrace{\left(\frac{1}{\lambda}\sum_{i=1}^{\ell-1}{\left(Y^i_t\right)}^{\T} \phi^i_t {\left(\phi_t^i\right)}^{\T} Y_t^i +  I_{d(n+1)} \right)}_{\overline{\Lambda}_t^{\ell}}.
\end{split}
\end{align*}
Then,
\begin{align}
\begin{split}
\frac{\det{\left(\Lambda_t^{\ell} \right)}}{\det{\left(\Lambda_t^{1} \right)}}&=\frac{\det{\left(\lambda \overline{\Lambda}_t^{\ell} \right)}}{\det{\left(\lambda I_{d(n+1)} \right)}}=\det{\left( \overline{\Lambda}_t^{\ell} \right)}\\
&=\det{\left( \frac{1}{\lambda}\sum_{i=1}^{\ell-1}{\left(Y^i_t\right)}^{\T} \phi^i_t {\left(\phi_t^i\right)}^{\T} Y_t^i +  I_{d(n+1)} \right)}\\
&=\prod_{i=1}^{d(n+1)} (1+\gamma_i),
\end{split}
\end{align}
where $\gamma_i$ is the $i$-th eigenvalue of $ \frac{1}{\lambda}\sum_{i=1}^{\ell-1}{\left(Y^i_t\right)}^{\T} \phi^i_t {\left(\phi_t^i\right)}^{\T} Y_t^i$. Let $\|{\left(Y^i_t\right)}^{\T} \phi^i_t \|\leq \delta_{\psi}$ for $t \in \{0,\cdots, T\}$ and $i\in\{1,\cdots,L\}$ (see \eqref{eq: bound on psi}). Then,
\begin{align}
\begin{split}
 \gamma_i &\leq \frac{1}{\lambda} \sum_{i=1}^{\ell-1}{\Tr{{\left(Y^i_t\right)}^{\T} \phi^i_t {\left(\phi_t^i\right)}^{\T} Y_t^i}}\\
  &\leq \frac{1}{\lambda} \sum_{i=1}^{\ell-1}\|{\left(Y^i_t\right)}^{\T} \phi^i_t\|^2\\
  &\leq \frac{\ell\delta_{\psi}^2}{\lambda}.
  \end{split}
\end{align}
Then
\begin{align}\label{eq: det_Lambda}
\begin{split}
\log{\left(\sqrt{ \det{\left( \overline{\Lambda}_t^{\ell} \right)}}\right)}&=\frac{1}{2}\log{\left( \prod_{i=1}^{d(n+1)} (1+\gamma_i)\right)}\\
&=\frac{1}{2} \sum_{i=1}^{d(n+1)} \log{\left(1+\gamma_i\right)}\\
&\leq \frac{d(n+1)}{2}\log{\left(1+\frac{\ell \delta_{\psi}^2}{\lambda}\right)}.
\end{split}
\end{align}
Substituting \eqref{eq: det_Lambda} in \eqref{eq: bounding r1}, we get with probability at least $\left.1-\delta\right.$
\begin{align}\label{eq: bounding r1 cont'd}
\begin{split}
&{\left\|\sum_{i=1}^{\ell-1}(\psi_t^i (y_{t}^i)^{\T} \left( v_{t+1}^{\ell}(s_{t+1}^i) - \Expect{}{ v_{t+1}^{\ell}(s_{t+1})}{s^i_t}\right)\right\|}^2_{\left({\Lambda_t^{\ell}}\right)^{-1}}\\
 &\qquad\leq  \sigma^2 \left(d(n+1)\log{\left(1+\frac{\ell \delta_{\psi}^2}{\lambda}\right)}+2\log{\left(\frac{1}{\delta} \right)}\right).
 \end{split}
\end{align}
Then,
\begin{align}\label{eq: bounding r1 final}
\begin{split}
&{\left\|\sum_{i=1}^{\ell-1}(\psi_t^i (y_{t}^i)^{\T} \left( v_{t+1}^{\ell}(s_{t+1}^i) - \Expect{}{ v_{t+1}^{\ell}(s_{t+1})}{s^i_t}\right)\right\|}_{\left({\Lambda_t^{\ell}}\right)^{-1}}\\
&\qquad \leq  \sigma \sqrt{d(n+1)\log{\left(1+\frac{\ell \delta_{\psi}^2}{\lambda}\right)}+2\log{\left(\frac{1}{\delta} \right)}}.
\end{split}
\end{align}
Then, we have
\begin{align}\label{eq: phi_Y_r1}
\begin{split}
 &\left|\phi_t^{\T}Y_t r_1\right| \\
 &\quad\leq \left\|\phi_t^{\T}Y_t (\Lambda_t^{\ell})^{-\frac{1}{2}} \right\| \\
 &\quad.{\left\|\sum_{i=1}^{\ell-1}(\psi_t^i (y_{t}^i)^{\T} \left( v_{t+1}^{\ell}(s_{t+1}^i) - \Expect{}{ v_{t+1}^{\ell}(s_{t+1})}{s^i_t}\right)\right\|}_{\left({\Lambda_t^{\ell}}\right)^{-1}}\\
 &\quad\leq
\sigma \sqrt{\left(d(n+1)\log{\left(1+\frac{\ell \delta_{\psi}^2}{\lambda}\right)}+2\log{\left(\frac{1}{\delta} \right)}\right)} \\
&\qquad  . \sqrt{\phi_t^{\T} Y_t \left(\Lambda_t^{\ell}\right)^{-1} Y_t^{\T} \phi_t}. 
 \end{split}
\end{align}
Next, we have
\begin{align}\label{eq: phi_Y_r3}
 \begin{split}
 \left|\phi_t^{\T}Y_t r_3\right| &\leq \lambda  \left\|\phi_t^{\T}Y_t \left(\Lambda_t^{\ell}\right)^{-\frac{1}{2}} \right\| \left\|\left(\Lambda_t^{\ell}\right)^{-\frac{1}{2}} \right\| \left\| \theta^*_{t+1} \right\|\\
 &\leq \sqrt{\lambda} \left\| \theta_{t+1}^* \right\|  \sqrt{\phi_t^{\T} Y_t \left(\Lambda_t^{\ell}\right)^{-1} Y_t^{\T} \phi_t}.
\end{split}
\end{align}
Next, we have
\begin{align}
\begin{split}
 r_2=&(\Lambda_t^{\ell})^{-1}  \Bigg(\sum_{i=1}^{\ell-1}(Y_t^i)^{\T} \phi_t^i(y_{t}^i)^{\T}\\
 &\qquad \quad .  \big( \Expect{}{v^{\ell}_{t+1}(s_{t+1})}{s^i_t}\!-\! \Expect{}{v^*_{t+1}(s_{t+1})}{s^i_t}\big)\Bigg)\\
 =&
 (\Lambda_t^{\ell})^{-1}  \Bigg(\sum_{i=1}^{\ell-1}(Y_t^i)^{\T} \phi_t^i (y_{t}^i)^{\T} \\
 &\quad .\int_{\mathcal{S}} \big(v^{\ell}_{t+1}(s_{t+1}) \!-\! v^{*}_{t+1}(s_{t+1})\big) \sum_{j=1}^{d} \phi^i_{j,t} \mu_{j,t}(ds_{t+1})\Bigg)\\
 =&
  (\Lambda_t^{\ell})^{-1}  \Bigg(\sum_{i=1}^{\ell-1}(Y_t^i)^{\T} \phi_t^i (y_{t}^i)^{\T} \\
  &\quad .\sum_{j=1}^{d} \phi_{j,t}^i \int_{\mathcal{S}} \left(v^{\ell}_{t+1}(s_{t+1}) \!-\! v^{*}_{t+1}(s_{t+1})\right)\mu_{j,t}(ds_{t+1})\Bigg)\\
  =&
    (\Lambda_t^{\ell})^{-1} \! \underbrace{\Big(\sum_{i=1}^{\ell-1}\psi_t^i  (\psi_t^i)^{\T} \Big)}_{ \Lambda_t^{\ell}-\lambda I_{d(n+1)}} \!\Expectover{\mu_{t+1}\!\!\!}{v^{\ell}_{t+1}(s_{t+1})\! -\! v^{*}_{t+1}(s_{t+1})}\\
   =&
   \Expectover{\mu_{t+1}}{v^{\ell}_{t+1}(s_{t+1})\! -\! v^{*}_{t+1}(s_{t+1})}\\
   &- \lambda  (\Lambda_t^{\ell})^{-1} \Expectover{\mu_{t+1}}{v^{\ell}_{t+1}(s_{t+1}) - v^*_{t+1}(s_{t+1})}\!.
 \end{split}
\end{align}
Then, we have
\begin{align}\label{eq: r2 final}
 \begin{split}
\phi_t^{\T}Y_t r_2\! =& \phi_t^{\T}Y_t  \Expectover{\mu_{t+1}}{v^{\ell}_{t+1}(s_{t+1}) - v^*_{t+1}(s_{t+1})}\\
 &\!-  \lambda  \phi_t^{\T}Y_t(\Lambda_t^{\ell})^{-1} \Expectover{\mu_{t+1}}{v^{\ell}_{t+1}(s_{t+1}) \!-\! v^*_{t+1}(s_{t+1})}\\
\!=& y_t^{\T}\Expectover{}{v^{\ell}_{t+1}(s_{t+1}) - v^*_{t+1}(s_{t+1}) | s_t}\\
\! &\!-  \lambda  \phi_t^{\T}Y_t(\Lambda_t^{\ell})^{-1} \Expectover{\mu_{t+1}}{v^{\ell}_{t+1}(s_{t+1}) \!-\! v^*_{t+1}(s_{t+1})}\\
=& \Expectover{}{V^{\ell}_{t+1}(x_{t+1},s_{t+1}) - V_{t+1}(x_{t+1},s_{t+1}) | s_t} \\
&\!\!-  \underbrace{\lambda  \phi_t^{\T}Y_t(\Lambda_t^{\ell})^{-1} \Expectover{\mu_{t+1}}{v^{\ell}_{t+1}(s_{t+1})\! -\! v^*_{t+1}(s_{t+1})}}_{r_4}.
\end{split}
\end{align}
Since we choose $R_{\theta}$ in Algorithm \ref{alg: alg1} to be the upper bound on $\|\theta_{t+1}^*\|$ (which we derive in Appendix \ref{subsec: R_theta}) for $t\in\{0,\!\cdots\!, T\!-\!1\}$, and we have $\|\theta^{\ell}_{t+1}\|\!\leq\! R_{\theta}$ (from Algorithm \ref{alg: alg1}) for $t\!\in\!\{0,\cdots, T-1\}$ and $\ell \!\in\! \{1,\!\cdots\! , L\}$, it can be seen from \eqref{eq: h bound} and \eqref{eq: q bound} that $\|v^{\ell}_{t+1}(s_{t+1})\|\leq\delta_{v}$ and $\|v^{*}_{t+1}(s_{t+1})\|\leq\delta_{v}$. We can derive an expression for $\delta_v$ using \eqref{eq: h bound} and \eqref{eq: q bound} as
\begin{align}
\begin{split}
\| v_{t+1}^{\ell}(s_{t+1}&)\| =\left\| 
\begin{bmatrix}
 h_{t+1}^{\ell}(s_{t+1})\\
 q_{t+1}^{\ell}(s_{t+1})
\end{bmatrix}\right\|\\
&\leq \underbrace{\sqrt{\|h_{t+1}^{\ell}(s_{t+1})\|^2 + \|q_{t+1}^{\ell}(s_{t+1})\|^2}}_{\delta_v}.
\end{split}
\end{align}
Then, we can bound $|r_4|$ in \eqref{eq: r2 final} as
\begin{align}\label{eq: r4 bound}
 \begin{split}
&\left| \lambda  \phi_t^{\T}Y_t(\Lambda_t^{\ell})^{-1} \Expectover{\mu_{t+1}}{v^{\ell}_{t+1}(s_{t+1}) - v^*_{t+1}(s_{t+1})} \right|\\
 &\quad\leq \lambda \left\| \phi_t^{\T}Y_t(\Lambda_t^{\ell})^{-\frac{1}{2}}\right\|  \left\|(\Lambda_t^{\ell})^{-\frac{1}{2}} \right\|\! \left\|  \Expectover{\mu_{t+1}}{v^{\ell}_{t+1} - v^*_{t+1}} \right\|\\
&\quad\leq 2 \sqrt{\lambda} \delta_{v} \sqrt{  \phi_t^{\T} Y_t\left(\Lambda_t^{\ell}\right)^{-1}Y_t^{\T} \phi_t}.
\end{split}
\end{align}
Using the parametrized form of the $Q$-function from Theorem \ref{thrm: Qfunction}, we have
\begin{align}\label{eq: delta Q}
\begin{split}
& Q^{\ell}_t(x_t,s_t,u_t)-  Q_t^*(x_t,s_t,u_t) = \phi^{\T}_t Y_t \left(\theta^{\ell}_{t+1} - \theta^*_{t+1}\right)\\
  &\qquad \qquad\qquad\qquad\qquad= \phi^{\T}_t Y_t \left(r_1 + r_2 + r_3\right)\\
&\implies Q^{\ell}_t(x_t,s_t,u_t)-  Q_t^*(x_t,s_t,u_t)\\
&\qquad \quad - \Expectover{}{V^{\ell}_{t+1}(x_{t+1},s_{t+1}) - V^*_{t+1}(x_{t+1},s_{t+1}) | s_t}\\
 &\qquad\qquad\qquad = \underbrace{\phi^{\T}_t Y_t \left(r_1 + r_3 \right) - r_4}_{\Delta_{t}^{\ell}(x_t,s_t,u_t)}.
 \end{split}
\end{align}
Then, using \eqref{eq: phi_Y_r1}, \eqref{eq: phi_Y_r3}, and \eqref{eq: r4 bound}, we can bound 
\begin{align}\label{eq: Delta bound}
 \begin{split}
 &\left|\Delta_{t}^{\ell}(x_t,s_t,u_t) \right| \\
 &\quad  \leq \Bigg(\sigma \sqrt{\left(d(n+1)\log{\left(1+\frac{\ell \delta_{\psi}^2}{\lambda}\right)}+2\log{\left(\frac{1}{\delta} \right)}\right)}\\
 &\qquad\qquad \!\!+ \left\| \theta^*_{t+1} \right\| \sqrt{\lambda}  + 2 \sqrt{\lambda} \delta_{v} \Bigg) \sqrt{\phi_t^{\T} Y_t \left(\Lambda_t^{\ell}\right)^{-1} Y_t^{\T} \phi_t}\\
 &\quad=\chi(\ell)\sqrt{\phi_t^{\T} Y_t \left(\Lambda_t^{\ell}\right)^{-1} Y_t^{\T} \phi_t}.
\end{split}
\end{align}
%
%
Let $\delta_t^{\ell} = V^{\ell}_t(x^{*\ell}_t, s^{\ell}_t) - V^*_t(x^{*\ell}_{t}, s^{\ell}_t) $ and $\zeta_{t+1}^{\ell} =\Expectover{}{\delta_{t+1}^{\ell}| s^{\ell}_t} - \delta_{t+1}^{\ell}$, where $x^{*\ell}_t$ is the state under the optimal policy $\pi^*_t$ starting from $x_0^{\ell}$ for $t\in \{0,\cdots, T\}$ and $\ell\in \{1,\cdots, L\}$. From the definition of the value function, we have $V_t^*(x,s)=\min_{u\in \mathcal{U}} Q_t^*(x,s,u)$. Further, since Algorithm \ref{alg: alg1} selects a greedy policy with respect to $Q_t^{\ell}(x,s,u)$, we have $V_t^{\ell}(x,s)=\min_{u\in \mathcal{U}} Q_t^{\ell}(x,s,u)$. Let $u_t^{*\ell}$ be the optimal control input at that generates $x_t^{*\ell}$ for $t\!\in\! \{0,\!\cdots\!, T\!-\!1\}$ and $\ell\in \{1,\cdots, L\}$. Then, we can write
\begin{align}\label{eq: delta_t}
\begin{split}
\delta_t^{\ell}&= V^{\ell}_t(x^{*\ell}_t, s^{\ell}_t) \!-\! V^*_t(x^{*\ell}_{t}, s^{\ell}_t) \\
 &= Q_t^{\ell}(x^{*\ell}_t, s^{\ell}_t,u_t^{\ell})\!-\!Q_t^{*}(x^{*\ell}_t, s^{\ell}_t,u_t^{*\ell})\\
 &\leq Q_t^{\ell}(x^{*\ell}_t, s^{\ell}_t,u_t^{*\ell})\!-\!Q_t^{*}(x^{*\ell}_t, s^{\ell}_t,u_t^{*\ell}).
 \end{split}
\end{align}
Then, from \eqref{eq: delta Q} we can write
\begin{align}\label{eq: Delta_Q at u_star}
\begin{split}
  &Q_t^{\ell}(x_t^{*\ell},s_t^{\ell},u^{*\ell})- Q_t^{*}(x_t^{*\ell},s_t^{\ell},u^{*\ell}) \\
  &\qquad \qquad \qquad= \Expectover{}{\delta^{\ell}_{t+1} | s_t^{\ell}}+\Delta_{t}^{\ell}(x_t^{*\ell},s_t^{\ell},u_t^{*\ell})\\
  &\qquad \qquad \qquad= \delta^{\ell}_{t+1} + \zeta_{t+1}^{\ell} + \Delta_{t}^{\ell}(x_t^{*\ell},s_t^{\ell},u_t^{*\ell}).
  \end{split}
\end{align}
Substituting  \eqref{eq: delta_t} in \eqref{eq: Delta_Q at u_star}, we get
\begin{align*}
 \delta_t^{\ell} \leq \delta^{\ell}_{t+1} + \zeta_{t+1}^{\ell}+ \Delta_{t}^{\ell}(x_t^{*\ell},s_t^{\ell},u_t^{*\ell}).
\end{align*}
Note that $x_0^{*\ell}=x_0^{\ell}$. Next, from \eqref{eq: regret} we write
\begin{align}\label{eq: regret bound proof}
\mathcal{R}(L)=\sum_{\ell=1}^L \delta_0^{\ell} \leq \sum_{\ell=1}^L \sum_{t=1}^T  \zeta_{t}^{\ell} + \sum_{\ell=1}^L \sum_{t=0}^{T-1}\Delta_{t}^{\ell}(x_t^{*\ell},s_t^{\ell},u_t^{*\ell}).
\end{align}
We have $\expect{\zeta_{t}^{\ell}|\mathcal{F}_{t-1}^{\ell}}=0$. Also, since we have $|\delta_t^{\ell}|\leq \sigma$ (see \eqref{eq: epsilon bound 2}) with $\sigma$ as in \eqref{eq: sigma bound}, then we have $|\zeta_t^{\ell}|\leq \sigma$ for $t\in\{0,\cdots,T\}$ and $\ell\in\{0,\cdots,L\}$. Hence $\zeta_{t}^{\ell}$ is a martingale difference sequence. Then, using the Azuma-Hoeffding inequality, for
any $\varepsilon>0$, we get
\begin{align*}
 \mathbb{P} \left( \sum_{\ell=1}^L \sum_{t=1}^T  \zeta_{t}^{\ell} \geq \varepsilon \right) \leq \exp\left(\frac{-\varepsilon^2}{2\sum_{i=1}^{LT} \sigma^2}\right) =\delta.
\end{align*}
Hence, we get with probability at least $1-\delta$
\begin{align}\label{eq: sum_zeta bound}
 \sum_{\ell=1}^L \sum_{t=1}^T  \zeta_{t}^{\ell} \leq \sigma \sqrt{2LT\log(1/\delta)}.
\end{align}
Next, using \eqref{eq: Delta bound} we can write
\begin{align*}
 &\sum_{\ell=1}^L \sum_{t=0}^{T-1}\Delta_{t}^{\ell}(x_t^{*\ell},s_t^{\ell},u_t^{*\ell})\\
 &\qquad \leq  \sum_{\ell=1}^L \sum_{t=0}^{T-1}\chi(\ell)\sqrt{(\phi_t^{\ell})^{\T} Y_t^{*\ell} \left(\Lambda_t^{\ell}\right)^{-1} (Y_t^{*\ell})^{\T} \phi_t(s_t^{\ell})}\\
 &\qquad \leq \chi(L) \sum_{\ell=1}^L \sum_{t=0}^{T-1}\|\left(\Lambda_t^{\ell}\right)^{-1/2}(Y_t^{*\ell})^{\T} \phi_t(s_t^{\ell})\|\\
 &\qquad \overset{\text(a)}{\leq} \chi(L) \delta_{\psi}\sum_{\ell=1}^L \sum_{t=0}^{T-1}\|\left(\Lambda_t^{\ell}\right)^{-1/2}\|,
\end{align*}
where in step (a) we have used \eqref{eq: bound on psi}. Assume the state transition matrix, $\varphi(t,0)$, in \eqref{eq: evolution of x} is nonsingular for $t\in\{0,\cdots,T\}$.\footnote{Since the system and the weight matrices are known, this assumption can be satisfied by an appropriate choice of the weight matrices.} Then, from Lemma \ref{lemma: expected_psi_psiT} we have $\expect{\psi_t \psi_t^{\T}}\succeq \gamma I_{d(n+1)}$ where $\gamma~=~\min_{t\in\{0,\cdots,T-1\}}\lambda_{\text{min}}\left(\expect{\psi_t\psi_t^{\T}}\right)~>~0$.
Further, from \eqref{eq: bound on psi}, we have $\|\psi_t^{\ell}\| \leq \delta_{\psi}$  for $t\in \{0,\cdots,T\}$ and $\ell \in \{1,\cdots,L\}$. Let $\delta \in [0,1]$ and assume $\ell \geq (8 \delta_{\psi}^2\log(d(n+1)/\delta))/\gamma$. Then, using Lemma \ref{lemma: norm_Lambda_inv} we have with probability at least $1-\delta$
\begin{align*}
\lambda_{\text{min}}(\Lambda_t^{\ell})\geq \lambda +\frac{(\ell-1)\gamma}{2}.
\end{align*}
Then, we can write 
\begin{align*}
  \sum_{\ell=1}^L \sum_{t=0}^{T-1} \|(\Lambda_t^{\ell}&)^{-1/2}\| \leq   \sum_{\ell=1}^L \sum_{t=0}^{T-1} \frac{1}{\sqrt{\lambda_{\text{min}}(\Lambda_t^{\ell})}}\\
  \leq&  \sum_{\ell=1}^L \frac{T}{\sqrt{\lambda +\frac{(\ell-1)\gamma}{2}}}\\
    \leq & \frac{T}{\sqrt{\lambda}} +\sum_{\ell=2}^L \frac{T}{\sqrt{\lambda +\frac{(\ell-1)\gamma}{2}}}\\
  \leq & \frac{T}{\sqrt{\lambda}} +\sum_{\ell=2}^L \frac{T\sqrt{2}}{\sqrt{(\ell-1)\gamma}}\\
    \leq & \frac{T}{\sqrt{\lambda}} + \frac{2\sqrt{2}T}{\sqrt{\gamma}}\sum_{\ell=2}^L( \sqrt{(\ell-1) - \sqrt{(\ell-2)})}\\
    = & \frac{T}{\sqrt{\lambda}} + \frac{2\sqrt{2}T}{\sqrt{\gamma}} \sqrt{(L-1)}\\
    \leq & \frac{T}{\sqrt{\lambda}} + \frac{4T\sqrt{L}}{\sqrt{\gamma}}.
\end{align*}
Then, we have 
\begin{align}\label{eq: Delta_star bound final}
 \sum_{\ell=1}^L \sum_{t=0}^{T-1}\Delta_{t}^{\ell}(x_t^{*\ell},s_t^{\ell},u_t^{*\ell})\leq \left(\frac{\delta_{\psi}T}{\sqrt{\lambda}} +\frac{4\delta_{\psi}T\sqrt{L}}{\sqrt{\gamma}} \right)\chi(L),
\end{align}
Substituting \eqref{eq: sum_zeta bound} and \eqref{eq: Delta_star bound final} in \eqref{eq: regret bound proof}, we get
\begin{align}\label{eq: regret bound proof final}
\mathcal{R}(L)\leq \sigma \sqrt{2LT\log(1/\delta)} \!+\! \left(\frac{\delta_{\psi}T}{\sqrt{\lambda}} \!+\!\frac{4\delta_{\psi}T\sqrt{L}}{\sqrt{\gamma}} \right)\chi(L).
\end{align}
the proof follows by substituting $\chi(L)$ defined in \eqref{eq: Delta bound} in \eqref{eq: regret bound proof final}, and the probability follows from the union bound.

\end{document}